\documentclass[11pt,a4paper]{article}
\usepackage[a4paper,margin=26mm,headheight=15pt]{geometry}
\usepackage[T1]{fontenc}
\usepackage{lmodern,microtype}
\DisableLigatures[-]{encoding=*,family=tt*}
\usepackage{amsmath,amssymb,amsthm,mathtools,mathrsfs}
\usepackage{booktabs,tabularx,array,longtable,multirow}
\usepackage{enumitem,algorithm,algpseudocode,float,caption,needspace}
\usepackage{hyperref,xcolor,fancyhdr}
\usepackage{tikz}
\usetikzlibrary{arrows.meta,positioning,fit,calc,decorations.pathreplacing}
\hypersetup{colorlinks=true,linkcolor=blue!45!black,citecolor=blue!45!black,urlcolor=blue!45!black,pdftitle={Recursive word conditioning for higher-dimensional partitions},pdfauthor={Hartosh Singh Bal}}
\allowdisplaybreaks[2]
\numberwithin{equation}{section}
\theoremstyle{plain}
\newtheorem{theorem}{Theorem}[section]
\newtheorem{proposition}[theorem]{Proposition}
\newtheorem{lemma}[theorem]{Lemma}
\newtheorem{corollary}[theorem]{Corollary}
\theoremstyle{definition}
\newtheorem{definition}[theorem]{Definition}

\theoremstyle{remark}
\newtheorem{remark}[theorem]{Remark}

\newcommand{\NN}{\mathbb N}
\newcommand{\wt}[1]{\lvert#1\rvert}
\newcommand{\subeq}{\preceq}
\newcommand{\supeq}{\succeq}
\newcommand{\Cont}{\operatorname{Cont}}
\newcommand{\eps}{\varepsilon}
\newcommand{\dd}{\mathrm d}
\title{\bfseries Recursive Word Conditioning for Higher-Dimensional Partitions:\\[3pt]Exact Enumeration and Complexity}
\author{Hartosh Singh Bal\\\small The Caravan, New Delhi}
\date{25 November 2026}
\begin{document}
\maketitle
\begin{abstract}
Let $p_d(n)$ denote the number of $d$-dimensional partitions of weight $n$. We develop an exact enumeration method operating on recursively parsed partition words (which are shown to encode $d$-dimensional partitions), beginning with an elementary interior-part decomposition of ordinary partitions. Conditioning on an interior word block separates the preceding and succeeding nested blocks. For higher-dimensional partitions, transversely re-parsing the preceding words permits the same conditioning operation to be iterated. Recursive containment and joins carry effective floors across the cuts, while a nonnegative excess permits exact truncation and pruning. The resulting procedure is a complete row-conditioning algorithm. For every fixed $d\ge3$, we prove that it computes $p_d(0),\ldots,p_d(N)$ using at most $\exp\!\bigl(O_d(N^{(d-1)/(d+1)}\log(N+2))\bigr)$ arithmetic operations. A diagonal refinement replaces the final row-conditioning stages by coupled interlacing chains of ordinary-partition words and sharpens the bound to $\exp\!\bigl(O_d(N^{(d-1)/(d+1)})\bigr)$. The algorithms count constrained completions through reusable weight-indexed coefficients rather than enumerating finished partitions. Exact computations with these algorithms reach weights 100, 50, 45, 40 and 35 in dimensions 3, 4, 5, 6 and 7, respectively. On a single 2.1\,GHz core the program computes $p_3(0),\ldots,p_3(100)$ in under four hours and each of the tables in dimensions 4--7 in under 30 minutes.   
\end{abstract}

\section{Introduction}\label{sec:intro}
The enumeration of higher-dimensional partitions remains difficult beyond the classical cases of ordinary and plane partitions. Let $p_d(n)$ denote the number of $d$-dimensional partitions of weight $n$. While Euler's and MacMahon's product formulas enumerate ordinary and plane partitions, respectively \cite{MacMahon}, MacMahon's proposed higher-dimensional extension does not give the correct coefficients. The rapid growth of admissible configurations makes exact counting challenging. For fixed $d$ the known estimates give $\log p_d(n)=\Theta_d(n^{d/(d+1)})$ \cite{BPA,Yeliussizov,Oganesyan}. Computing coefficients without generating every partition individually is therefore a natural objective.

The computational history of solid partitions illustrates this difficulty. Early methods include those of Atkin, Bratley, Macdonald and McKay, and Knuth \cite{AtkinEtAl,Knuth}. Mustonen and Rajesh extended exact enumeration through weight 50 \cite{MustonenRajesh}. Balakrishnan, Govindarajan and Prabhakar subsequently developed shape-equivalence methods and parallel computation, extending the sequence through weight 68; an IIT Madras computation reached weight 72 in January 2013 \cite{BGP,SPProject}. Other approaches, including Govindarajan's structural decompositions and the more recent socle-reduction method of Graffeo, Monavari, Moschetti and Ricolfi, advance enumeration through different reductions \cite{Govindarajan,GMMR2026}. The present work instead uses recursive conditioning of nested partition chains.

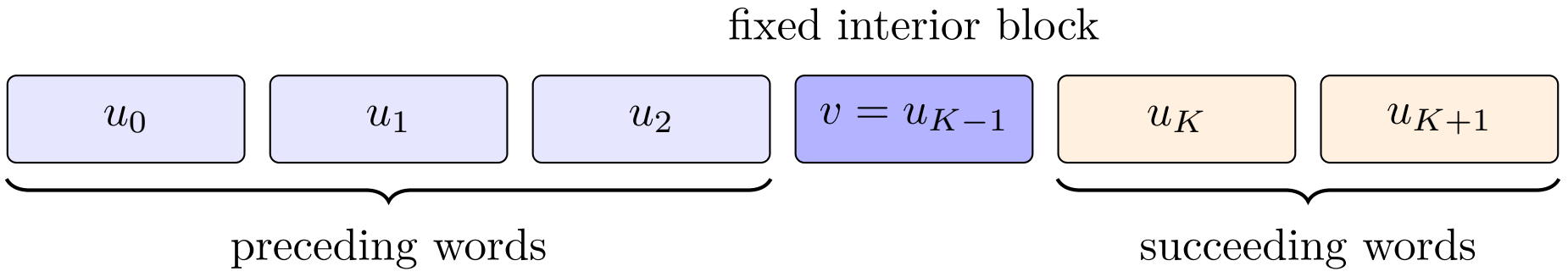
\begin{figure}[htbp]
\centering
\begin{tikzpicture}[>=Latex, font=\small, x=1cm, y=1cm,
  blk/.style={draw, rounded corners=2pt, minimum width=1.85cm, minimum height=.68cm, inner sep=2pt}]
\foreach \x/\lab/\nm in {0/{u_0}/b0, 1/{u_1}/b1, 2/{u_2}/b2}
  \node[blk, fill=blue!10] (\nm) at (2.05*\x,0) {$\lab$};
\node[blk, fill=blue!30] (bv) at (2.05*3,0) {$v=u_{K-1}$};
\node[blk, fill=orange!12] (b4) at (2.05*4,0) {$u_K$};
\node[blk, fill=orange!12] (b5) at (2.05*5,0) {$u_{K+1}$};
\node[above=4pt of bv] {fixed interior block};
\draw[decorate, decoration={brace, amplitude=5pt, mirror}, thick]
  ($(b0.south west)+(0,-.12)$) -- ($(b2.south east)+(0,-.12)$)
  node[midway, below=7pt] {preceding words};
\draw[decorate, decoration={brace, amplitude=5pt, mirror}, thick]
  ($(b4.south west)+(0,-.12)$) -- ($(b5.south east)+(0,-.12)$)
  node[midway, below=7pt] {succeeding words};
\end{tikzpicture}
\caption{An interior cut of a recursively parsed word (schematic). The preceding blocks contain the fixed block $v$ and force at least $K|v|$ total weight together with $v$; the succeeding blocks are bounded above by $v$. The displayed positions stand for arbitrary block indices, not a cut restricted to six blocks.}
\label{fig:interior-word-cut}
\end{figure}

The organizing operation in this paper is \emph{interior word conditioning}: instead of generating a partition word from its first block (Section~\ref{sec:words} shows that these words recursively encode $d$-dimensional partitions), fix a block in its unique recursive parse. The blocks preceding the cut must contain it and therefore carry compulsory weight; the succeeding blocks form a bounded continuation (Figure~\ref{fig:interior-word-cut}). The resulting exact enumeration is of \emph{residual freedom above a forced minimum}, not of completed partitions one by one. Re-reading the preceding blocks \emph{transversely} produces a decreasing chain of coupled lower-dimensional words on which the same conditioning operation can be repeated. The re-reading itself preserves all information; it is the next conditioning cut that bounds an additional index direction. Recursive joins encode inherited floors, and nonnegative excess weights make the remaining budget explicit. This yields a complete row-conditioning algorithm. A second reading uses diagonal blocks, again ordinary partition words, but with interlacing in place of containment. Its stronger compulsory completion weights lead to the improved operation bound. In both algorithms the states remain partition words: the row construction uses recursive containment between parsed blocks, while the diagonal construction uses interlacing between ordinary-partition words.

Our principal results are as follows.
\begin{theorem}[Recursive row conditioning]\label{thm:row-main}
For every fixed integer $d\ge3$, an explicitly specified deterministic row-conditioning algorithm computes $p_d(0),\ldots,p_d(N)$ exactly using at most
\begin{equation}\label{eq:row-headline}
\exp\!\left(O_d\!\left(N^{\frac{d-1}{d+1}}\log(N+2)\right)\right)
\end{equation}
arithmetic operations.
\end{theorem}
\begin{theorem}[Diagonal refinement]\label{thm:diag-main}
For every fixed integer $d\ge3$, replacing the final row-conditioning stages by a coupled diagonal transfer gives an exact algorithm computing the same coefficients using at most
\begin{equation}\label{eq:headline}
\exp\!\left(O_d\!\left(N^{\frac{d-1}{d+1}}\right)\right)
\end{equation}
arithmetic operations.
\end{theorem}
For solid partitions, Theorem~\ref{thm:diag-main} gives $\exp(O(\sqrt N))$. The diagonal refinement removes the logarithmic factor in the proved row-only bound; this does not establish that the factor is intrinsic to all row-only algorithms.

A direct transfer over successive slices has as its potential states the $(d-1)$-dimensional partitions of weight at most $N$, a collection of size $\exp(\Theta_d(N^{(d-1)/d}))$. Interior conditioning reduces the exponent from $(d-1)/d$ to $(d-1)/(d+1)$, up to a logarithmic factor for Algorithm R and without it for Algorithm D. Since $\log p_d(N)=\Theta_d(N^{d/(d+1)})$, the diagonal algorithm computes the entire coefficient table through $N$ using
\[
 p_d(N)^{\,O_d(N^{-1/(d+1)})}
\]
arithmetic operations, subpolynomial in its final coefficient.

Section~\ref{sec:computation} lists the coefficients obtained with these algorithms, through $p_3(100)$, $p_4(50)$, $p_5(45)$, $p_6(40)$ and $p_7(35)$, and compares them with the data of Graffeo, Monavari, Moschetti and Ricolfi \cite{GMMR2026,GMMRData}.

The method has several antecedents but differs in how its counting problems are organised. Classical Durfee decompositions use distinguished interior structures \cite{AndrewsDurfee}; transfer-matrix and umbral methods maintain boundary data across configurations \cite{Zeilberger}; and diagonal interlacing is classical in plane partitions and the Schur process \cite{OkounkovReshetikhin}. Related enumerations with prescribed shapes or boundaries appear in the work of Stanley, Gansner, Sagan and Morales--Pak--Panova, and in topological-vertex constructions \cite{StanleyOrdered,Gansner,SaganSkew,MPP,TopVertex}. Amanov and Yeliussizov give a MacMahon-type product for higher-dimensional partitions counted by a corner-hook statistic rather than by volume \cite{AmanovYeliussizov}. Earlier solid-partition computations exploited sophisticated shape grouping, while recent all-dimensional methods use embedding-dimension and socle reduction. The contribution here is the systematic combination of interior conditioning, effective floors, recursive dimension reduction and reusable constrained completion counts, together with uniform complexity bounds. In particular, the method addresses completion above arbitrary finite floors rather than relying on an unrestricted product formula.

The work is rooted in the author's earlier recursive partition-word construction \cite{BalWords}. Although the construction can now be re-expressed in standard language for $d$-dimensional partitions, we retain the word formulation that led to the algorithm. Its fully realised enumerative consequences were developed through trial-and-error coding experiments, algorithmic design, proof development and verification undertaken in sustained human--AI research interactions. The specific roles of the AI research assistants are disclosed in the acknowledgments.

Section~\ref{sec:words} introduces recursive words, their operations and the ordinary-partition prototype; Section~\ref{sec:example} motivates the interior cut through a hand calculation. Sections~\ref{sec:conditioning} and \ref{sec:alg} establish the floor-constrained recurrence and complete row algorithm. Section~\ref{sec:diagonal} develops the diagonal refinement; Section~\ref{sec:complexity} proves the bounds. Section~\ref{sec:discussion} discusses the scope of the results, and Section~\ref{sec:computation} lists the computed coefficients; the appendices follow, the last of them describing the ancillary files.

\section{Recursive partition words and nested chains}\label{sec:words}
Before the general grammar, consider two examples. The ordinary partition $\lambda=(3,2)$ has word $1101$: each run of ones records one less than the corresponding part, and a zero separates parts. The plane partition with rows $(2,1)$ and $(1)$ has word $210$: its first row has binary word $10$, shifted to $21$, followed by a zero and the shifted empty word of its second row. In particular, the \emph{empty word} denotes one cell; it must not be confused with the null partition.

Throughout, $\NN=\{0,1,2,\ldots\}$ and $\wt{u}$ denotes word weight. We retain the convention that an ordinary, plane or solid partition has respectively one, two or three independent part indices.

\begin{definition}[Recursive partition words]\label{def:words}
Let $\varnothing$ be a distinguished \emph{null word} of weight zero, and let $\epsilon$ be the empty string, representing the unique one-cell partition and having weight one. For an ordinary partition $\lambda=(\lambda_1\ge\cdots\ge\lambda_r>0)$, its legal binary word is
\[
 W_1(\lambda)=1^{\lambda_1-1}0\,1^{\lambda_2-1}0\cdots0\,1^{\lambda_r-1}.
\]
Equivalently, the successive runs of ones, including empty runs, have lengths $\lambda_i-1$ for a weakly decreasing sequence of positive integers $\lambda_i$.

Inductively, a nonnull legal $k$-word for $k\ge2$ has the unique parse
\[
 w=\iota(u_0)\,0\,\iota(u_1)\,0\cdots0\,\iota(u_{r-1}),
 \qquad u_0\succeq\cdots\succeq u_{r-1}\ne\varnothing,
\]
where the $u_i$ are legal $(k-1)$-words, the order $\succeq$ is defined recursively below, and $\iota$ increases every letter by one (fixing $\epsilon$). The outer zeroes separate shifted blocks; each intervening block may be empty as a string but encodes a nonnull one-cell partition. A legal nonnull word of length $n-1$ has weight $\wt w=n$, so $\wt w=\operatorname{length}(w)+1$. We write $W_k(\pi)$ for the canonical $k$-word encoding a $k$-partition $\pi$ (and retain $W_k(\varnothing)=\varnothing$ for the null partition); in particular $W_k$ denotes this same recursive block assembly at every dimension.
\end{definition}

\begin{definition}[Recursive containment and join]\label{def:word-operations}
The null word is contained in every legal word, and $\varnothing\vee u=u$. For ordinary binary words, parse the part sequences $\lambda=(\lambda_i)$ and $\mu=(\mu_i)$, padding missing parts by zeros. Define $u\supeq v$ when $\lambda_i\ge\mu_i$ for every $i$, and define $u\vee v$ by the partition with parts $\max(\lambda_i,\mu_i)$. For $k\ge2$, parse nonnull words into lower-dimensional blocks $u_i$ and $v_i$, padding absent blocks by $\varnothing$, and put
\[
 u\supeq v\quad\Longleftrightarrow\quad
 u_i\supeq v_i\ \text{for every }i,
 \qquad
 u\vee v=\operatorname{assemble}_k(u_0\vee v_0,u_1\vee v_1,\ldots).
\]
Here $\operatorname{assemble}_k$ shifts and joins the consecutive nonnull blocks with outer zero separators, omitting the trailing null blocks. For an ordered finite family, its weight is the sum of the word weights of its members.
\end{definition}

\begin{proposition}[Closure of the word operations]\label{prop:wordops}
Recursive containment is a partial order on legal $k$-words. Recursive join is legal and is the least word containing both arguments. For $k\ge2$, the weight of any nonnull $k$-word is the sum of the weights of its parsed lower-dimensional blocks; for $k=1$, it is the sum of its parts. Thus containment, join and weight can all be evaluated on recursively parsed words, without constructing a higher-dimensional part array.
\end{proposition}
\begin{proof}
For ordinary partitions the claims are immediate from the coordinatewise order on their padded part lists: taking componentwise maxima preserves weak decrease. Induct on $k$. The blockwise joins of two decreasing block sequences are again decreasing, since join is monotone in each argument. After the trailing null blocks are dropped, their assembly is therefore a legal word, and blockwise minimality makes it the least upper bound. Finally, each lower block of weight $m$ contributes $m-1$ letters, and the $r$ nonnull blocks contribute $r-1$ separators; their assembled word has length $\sum_{i=0}^{r-1}(m_i-1)+(r-1)=\sum_i m_i-1$. This proves the weight identity.
\end{proof}

\begin{definition}[Equivalent part-array interpretation]\label{def:kpartition}
A $k$-partition may equivalently be regarded as a finitely supported function $\pi:\NN^k\to\NN$ weakly decreasing in each coordinate, with weight $\wt\pi=\sum_{x\in\NN^k}\pi(x)$. The $k$-word is the recursive encoding of its successive nonzero coordinate slices $\pi_0\supeq\pi_1\supeq\cdots$; its word containment, join and weight agree with pointwise containment, pointwise maximum and total volume. This correspondence follows inductively by matching zero-separated blocks with coordinate slices, starting from the padded part lists in dimension one. The corresponding finite order ideal is
\[
 I(\pi)=\{(x,h)\in\NN^{k+1}:0\le h<\pi(x)\}.
\]
We write $p_k(n)$ for the number of words, or equivalently part arrays, of weight $n$, with $p_k(0)=1$. We use arrays and ideals when a geometric argument or an established counting bound is more convenient.
\end{definition}

The recursive grammar gives an intrinsic legality test: each shifted block is a legal lower-dimensional word, and consecutive blocks are nested. The weight and order conventions remain valid when one or several trailing blocks are absent, because these are padded by $\varnothing$, never by $\epsilon$. At $k=3$, reading from left to right, the letter $3$ adds a cell to the current pile, $2$ starts a new pile in the current row, $1$ starts a new row of the current slice, and $0$ starts a new slice, as in the word $321020$ of Appendix~\ref{app:small}; every intermediate word must be legal.

\begin{proposition}[Transverse re-parsing of coupled words]\label{prop:transverse}
Let $k\ge2$, let $\Gamma$ be a finite product of ordered chains, and let $(u_\gamma)_{\gamma\in\Gamma}$ be a coupled family of legal $k$-words, meaning $u_\gamma\supeq u_{\gamma'}$ whenever $\gamma\le\gamma'$. Parse each word into its $(k-1)$-blocks $u_{\gamma,r}$, padded by null blocks at absent positions. The transverse families
\[
 c_r=(u_{\gamma,r})_{\gamma\in\Gamma}\qquad(r\ge0)
\]
are coupled families of legal $(k-1)$-words, and
\[
 c_0\supeq c_1\supeq\cdots,\qquad
 \sum_{\gamma\in\Gamma}\wt{u_\gamma}
 =\sum_{r\ge0}\sum_{\gamma\in\Gamma}\wt{u_{\gamma,r}}.
\]
Conversely, every eventually-null decreasing chain $(c_r)$ of coupled $(k-1)$-word families uniquely reassembles into a coupled $k$-word family of the same total weight. When $\Gamma$ is a single chain of length $J$, each $c_r$ itself encodes a legal $k$-word with at most $J$ nonnull outer blocks; for a general product grid it is retained as a coupled family of words.
\end{proposition}
\begin{proof}
The two inequalities defining the coupled family and its nested blocks are exactly
$u_{\gamma,r}\supeq u_{\gamma',r}$ for $\gamma\le\gamma'$ and
$u_{\gamma,r}\supeq u_{\gamma,r+1}$. These become, respectively, coupling inside each $c_r$ and decrease from $c_r$ to $c_{r+1}$. Conversely the same comparisons guarantee legality and coupling of the reassembled words. Unique zero-separator parsing gives a bijection; the weight identity is Proposition~\ref{prop:wordops} summed over the family. The single-chain statement follows by assembling its nested blocks in that order.
\end{proof}

\begin{remark}[Canonical word states]\label{rem:wordtree}
Conditioning at an interior block and conditioning again on a transverse block are thus the same operation on recursively nested words. Compact word addresses can record the blocks, floors and joins; they need not be literal character strings. The interval generator in Section~\ref{sec:alg} gives the required duplicate-free legal extensions.
\end{remark}

\subsection{The ordinary-partition prototype: conditioning at the \texorpdfstring{$K$}{K}th part}\label{subsec:ordinary-prototype}
The principle behind the higher-dimensional algorithm is already exact for binary partition words. Fix $K\ge2$ and parse the canonical word of an ordinary partition into its successive runs of ones, separated by zeroes. A run of length $a-1$ encodes a part of size $a$. Padding the part list by zeros, fix the $K$th part to be $a\ge0$. If $a>0$, the $K-1$ preceding runs encode parts at least $a$, and the succeeding runs encode parts at most $a$. Thus a word with $\lambda_K=a>0$ is uniquely determined by three pieces: an arbitrary partition $\mu$ with at most $K-1$ parts, obtained by subtracting $a$ from each of the first $K-1$ parts; the fixed $K$th part $a$; and an arbitrary partition $\nu$ whose largest part is at most $a$. Conversely, adding $a$ to the padded parts of $\mu$, inserting the fixed part, and adjoining $\nu$ uniquely reassembles a legal binary word. Its weight is
\begin{equation}\label{eq:ordinary-weight-split}
 |\lambda|=Ka+|\mu|+|\nu|.
\end{equation}
The first term is compulsory; the other two terms count independent residual choices. The familiar finite Euler products therefore give a complete elementary conditioning identity.
\begin{proposition}[Interior-part identity]\label{prop:ordinary-interior}
For every integer $K\ge2$, with $(q;q)_m=\prod_{j=1}^{m}(1-q^j)$ and $(q;q)_0=1$,
\begin{equation}\label{eq:ordinary-interior}
 \prod_{j\ge1}\frac1{1-q^j}
 =\frac1{(q;q)_{K-1}}
 +\sum_{a\ge1}\frac{q^{Ka}}{(q;q)_{K-1}(q;q)_a}.
\end{equation}
\end{proposition}
\begin{proof}
The partitions with fewer than $K$ nonzero parts are counted by $(q;q)_{K-1}^{-1}$, by conjugation. Every other partition has a unique positive $K$th part $a$. The preceding residual partition $\mu$ has at most $K-1$ parts, and the suffix $\nu$ has largest part at most $a$; their generating functions are $(q;q)_{K-1}^{-1}$ and $(q;q)_a^{-1}$, respectively. Equation~\eqref{eq:ordinary-weight-split} supplies the factor $q^{Ka}$. The cases are disjoint and exhaustive, proving the identity.
\end{proof}
Identity \eqref{eq:ordinary-interior} is Euler's identity $\sum_{a\ge0}q^{Ka}/(q;q)_a=\prod_{j\ge K}(1-q^j)^{-1}$ \cite[Ch.~2]{AndrewsBook} multiplied by $1/(q;q)_{K-1}$; the proof above reads it combinatorially at the $K$th part.

The cut position $K$ is mathematically arbitrary: changing it changes the constrained states but not the generating function. In higher dimensions the fixed part $a$ is replaced by a recursively parsed word block $v$; the preceding and succeeding words are again independently reassembled once $v$ is fixed. Unlike this elementary case, however, their constrained completion series above arbitrary floors need not factor into Euler products. This is why the general algorithm stores exact completion coefficients rather than invoking an unrestricted product identity.

\section{Why an interior cut helps: fifty-four configurations}\label{sec:example}
The ordinary-partition identity of Proposition~\ref{prop:ordinary-interior} gives the first exact interior cut. We now see the same idea at work in a solid-partition calculation small enough to check directly, before introducing the general floor-constrained formalism. A solid partition is a decreasing chain of plane partitions. Fix its third slice to be the one-cell plane partition $v=(1)$ and consider the case with no subsequent nonempty slices. At total weight seven, we must count
\begin{equation}\label{eq:54-start}
A\succeq B\succeq v,\qquad |A|+|B|=6.
\end{equation}
The compulsory three slices $(v,v,v)$ use three cells, leaving an excess budget of four. We count the same configurations in two ways.

\paragraph{Condition on the complete intermediate plane partition.}
Choose $B$ first. Since $A\succeq B\succeq v$, its weight is 1, 2 or 3. If $|B|=1$, then $B=v$ and any of the $p_2(5)=24$ plane partitions of weight five can serve as $A$. If $|B|=2$, there are three shapes for $B$, each contained in exactly eight plane partitions $A$ of weight four. If $|B|=3$, equal weights force $A=B$, giving $p_2(3)=6$ choices. Hence
\begin{equation}\label{eq:54-whole}
24+3\cdot8+6=54.
\end{equation}
This procedure conditions on ten complete intermediate plane partitions: one of weight one, three of weight two, and six of weight three.

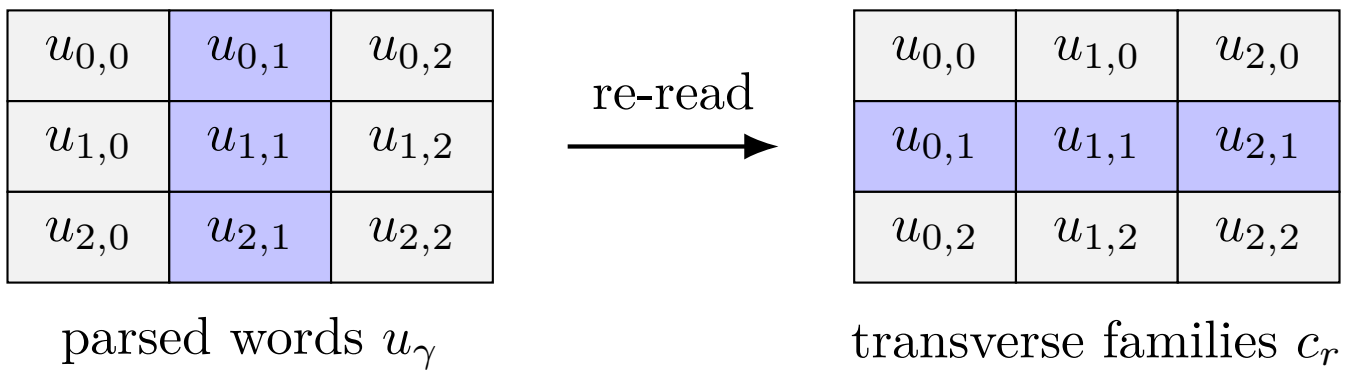
\begin{figure}[htbp]
\centering
\begin{tikzpicture}[>=Latex,font=\small,
  cell/.style={draw,minimum width=1.07cm,minimum height=.6cm,inner sep=0pt,outer sep=0pt}]
\foreach \i in {0,1,2}{
  \foreach \j in {0,1,2}{
    \ifnum\j=1 \def\shade{blue!23}\else\def\shade{gray!10}\fi
    \node[cell,fill=\shade] at (1.07*\j,-.6*\i) {$u_{\i,\j}$};
  }
}
\draw[->,thick] (3.17,-.6)--(4.57,-.6) node[midway,above=3pt]{re-read};
\foreach \j in {0,1,2}{
  \foreach \i in {0,1,2}{
    \ifnum\j=1 \def\shade{blue!23}\else\def\shade{gray!10}\fi
    \node[cell,fill=\shade] at (5.6+1.07*\i,-.6*\j) {$u_{\i,\j}$};
  }
}
\node at (1.07,-1.9) {parsed words $u_\gamma$};
\node at (6.67,-1.9) {transverse families $c_r$};
\node[align=center,font=\footnotesize] at (3.87,-2.5) {The same blocks, re-indexed; a later cut, not the re-reading, bounds an index direction.};
\end{tikzpicture}
\caption{Transverse re-parsing of three coupled words with three displayed blocks each. The shaded column becomes the shaded transverse family. The bijection preserves the blocks, their weight, and both directions of coupling.}
\label{fig:transverse-reparse}
\end{figure}

\paragraph{Condition again, now on a transverse word.}
Parse the two plane-partition words $A$ and $B$ into their ordinary-partition word blocks $A=(A_0,A_1,\ldots)$ and $B=(B_0,B_1,\ldots)$. Re-read corresponding blocks \emph{across} the two words rather than along either word separately (Figure~\ref{fig:transverse-reparse}):
\[
 c_r=(A_r,B_r),\qquad c_0\supeq c_1\supeq\cdots.
\]
Because $A_r\supeq B_r$, every $c_r$ is itself the block sequence of a legal plane-partition word with at most two nonnull outer blocks; $c_r=(\varnothing,\varnothing)$ represents the null word. This is the transverse re-parsing of Proposition~\ref{prop:transverse}. Fix its \emph{second transverse block} $c=c_1=(A_1,B_1)$, exactly as the first cut fixed an interior block of the original word. In the table below, partition shapes label the constituent ordinary-partition words, while $W_2(c)$ displays their assembled plane-partition word. In particular, $(1)$ denotes $\epsilon$ as an ordinary-partition word, whereas $\varnothing$ denotes a missing block.

For each $c$, the earlier transverse word $c_0=(A_0,B_0)$ and the later transverse words $c_r=(A_r,B_r)$ for $r\ge2$ can be counted separately. The earlier rows must contain both $c$ and the first-row floor inherited from $v$. Their least possible values are therefore componentwise joins. The weight forced by this choice is its \emph{compulsory excess} $\varepsilon(c)$; it is computed relative to the baseline pair $(v,v)$. All admissible $c$ for \eqref{eq:54-start}, with their exact contributions, are as follows:
\begin{center}\small
\begin{tabular}{@{}llrr@{}}\toprule
Fixed transverse pair $c$ & Assembled word $W_2(c)$ & $\varepsilon(c)$ & Contribution\\\midrule
$(\varnothing,\varnothing)$ & $\varnothing$ & 0 & 18\\
$((1),\varnothing)$ & $\epsilon$ & 1 & 17\\
$((1),(1))$ & $0$ & 2 & 9\\
$((2),\varnothing)$ & $2$ & 3 & 4\\
$((1,1),\varnothing)$ & $1$ & 3 & 4\\
$((2),(1))$ & $20$ & 4 & 1\\
$((1,1),(1))$ & $10$ & 4 & 1\\\midrule
\multicolumn{3}{r}{Total} & \textbf{54}\\\bottomrule
\end{tabular}
\end{center}
We have replaced ten complete plane-partition conditioning states by seven bounded transverse plane-partition words (equivalently, coupled ordinary-word pairs). The small numerical reduction is not itself a speedup argument; the gain arises from the bounded transverse extent, forced weight bounds and reuse of completion counts.

\paragraph{Deriving the contribution seventeen.}
For $c=((1),\varnothing)$, the fixed second row of $A$ uses one excess cell. There are three cells left. The first rows satisfy
\[
A_0\succeq B_0\succeq (1).
\]
Below the cut, $B$ has no nonzero rows, while $A$ may have any number of further one-cell rows. Thus each unused cell from the first-row pair determines a unique lower continuation. Fix $B_0$ and count the possibilities for $A_0$ by its weight.

If $B_0=(1)$, the admissible weights of $A_0$ are 1, 2, 3 and 4. The numbers of ordinary partitions at those weights are 1, 2, 3 and 5, respectively; the remaining budget is spent uniquely below the cut. The contribution is therefore $1+2+3+5=11$. If $B_0=(2)$, then $A_0$ is either $(2)$ at weight two or $(3),(2,1)$ at weight three, giving three possibilities. For $B_0=(1,1)$, the possible $A_0$ are $(1,1)$ at weight two and $(2,1),(1,1,1)$ at weight three, again giving three. No larger $B_0$ fits the budget. Hence
\begin{equation}\label{eq:17-byweight}
11+3+3=17.
\end{equation}
Equivalently, group the first-row pairs by their excess $e=0,1,2,3$ above the compulsory two one-cell rows. Their counts are
\[
U_c=(1,2,5,9).
\]
The lower-chain counts are $L_c=(1,1,1,1)$, so their convolution at residual weight three is $1+2+5+9=17$. These are \emph{stored coefficient vectors}, not lists of completed plane partitions. The method is direct integer counting; generating-function products merely abbreviate the convolution.

\paragraph{The other six states.}
For $c=(\varnothing,\varnothing)$ there are no later rows, and the first-row pair of total weight six contributes $7+8+3=18$, grouped by $|B_0|=1,2,3$. If $c=((1),(1))$, its compulsory cost is two, leaving excess two; the upper vector $(1,2,5)$ and lower vector $(1,1,2)$ give $5+2+2=9$. For either $c=((2),\varnothing)$ or $c=((1,1),\varnothing)$, the conditioned row and forced first-row enlargement cost three, leaving one cell: three upper enlargements and one lower continuation, contributing four. The remaining two conditioned rows already use all four excess cells and contribute one each. Thus the second decomposition gives
\begin{equation}\label{eq:54-rows}
18+17+9+4+4+1+1=54.
\end{equation}

The two principles visible here scale beyond the example. First, the original cut and this transverse cut are applications of \emph{one word-conditioning rule} to differently read nested word families. At the second cut, the transverse direction has length at most $K-1$, while the dimension of each constituent block is one lower. Repeating the re-parsing produces product-grid-indexed coupled words, with another bounded direction at each cut. Second, if the transverse cut is at depth $S$, the preceding $S-1$ transverse words must contain $c$; together with $c$ itself this forces $S\wt{c}\le N$. This \emph{forced-smallness} restriction accompanies the new bounded direction. For a fixed number of members, ordinary-partition states of weight at most $m$ have counting-logarithm scale $O(m^{1/2})$, compared with $O(m^{2/3})$ for complete plane partitions. These state-population comparisons motivate the method but do not by themselves prove its total operation bound.

\paragraph{A nontrivial effective floor.}
To see why the join cannot be suppressed, take an inherited plane-partition floor with first two rows $v_0=(2)$ and $v_1=(1)$. Fix the second-row pair $c=((1,1),(1))$. The least permissible first rows are
\[
A_0^{\min}=(1,1)\vee(2)=(2,1),\qquad
B_0^{\min}=(1)\vee(2)=(2).
\]
The least pair containing the floor has total weight six, whereas the least pair with the cut fixed has weight eight. The compulsory excess is therefore two. Any valid count at total budget $E$ must first pay this cost, then split the remaining $E-2$ between upper and lower completions. The next section formalises this accounting.

The complete weight-seven coefficient $p_3(7)=307$, including third-slice choices outside the single-cell branch, is reconstructed in Appendix~\ref{app:small}.

\section{Floor-constrained families and the interior-cut recurrence}\label{sec:conditioning}
\begin{definition}[Floored coupled word family]\label{def:family}
Let $\Gamma=[g_1]\times\cdots\times[g_t]$, where $[g]=\{0,\ldots,g-1\}$ and $t=0$ is allowed. A $\Gamma$-family of legal $k$-words is $u=(u_\gamma)_{\gamma\in\Gamma}$ with $u_\gamma\supeq u_{\gamma'}$ whenever $\gamma\le\gamma'$, using recursive word containment from Definition~\ref{def:word-operations}. Put $\wt{u}=\sum_\gamma\wt{u_\gamma}$. Given a coupled word family $\phi$, called a \emph{floor}, define
\begin{equation}\label{eq:Cdef}
 C_k(\Gamma,\phi;q)=\sum_{\substack{u\text{ a }\Gamma\text{-family}\\u\supeq\phi}}q^{\wt{u}-\wt{\phi}}.
\end{equation}
The constant term is one. The family may equivalently be described as $\wt{\Gamma}$ coupled $k$-partitions, or as a finite order ideal in $\Gamma\times\NN^{k+1}$ containing the floor ideal. The counting states and all order tests, however, are already defined in the recursive word language.
\end{definition}

The equivalent product-order description has $t+k+1$ coordinate directions. Transverse re-parsing first reorganises the unbounded block index into a chain of coupled words. It is the subsequent interior cut, retaining its $S-1$ preceding transverse blocks, that turns this portion of the chain into an additional bounded index-grid direction while lowering the dimension of its constituent words by one. In particular, the conditioned $(k-1)$-word families on a grid with $t$ axes have a product-order description with $t+k=d$ directions at the relevant state-enumeration level. This fixed-dimensional identity will be used for the state-counting proof; it does not require the algorithm to construct an ideal.

\begin{definition}[Floored lower chains]\label{def:chain}
Let $\Phi=(\Phi_r)_{r\ge r_0}$ be a weakly decreasing, eventually-null sequence of coupled $\Gamma$-families of legal $k$-words. For $X\supeq\Phi_r$ define
\begin{equation}\label{eq:Fdef}
 F(X,r;\Phi;q)=\sum_{\substack{X\supeq Y_r\supeq Y_{r+1}\supeq\cdots\\Y_j\supeq\Phi_j\text{ for every }j}}
 q^{\sum_{j\ge r}\wt{Y_j}},\qquad
 G(X,r;\Phi;q)=q^{-\sum_{j\ge r}\wt{\Phi_j}}F(X,r;\Phi;q).
\end{equation}
All chains are eventually zero. Every exponent in $G$ is nonnegative. If all floors vanish from $r$ onward, we abbreviate $F(X,r;0;q)$ to $F(X;q)$.
\end{definition}

\subsection{Conditioning on a block of the original word}
Let $P_d(q)=\sum_{n\ge0}p_d(n)q^n$. For $K\ge2$ put $J=K-1$ and let $v^{(J)}$ denote the $[J]$-family all of whose members equal $v$.

\begin{lemma}[Block-conditioning identity]\label{lem:top}
For every $d\ge2$ and $K\ge2$,
\begin{equation}\label{eq:top}
 P_d(q)=\sum_{v}q^{K\wt{v}}\,C_{d-1}([J],v^{(J)};q)\,F(v;q),
\end{equation}
where the sum ranges over $(d-1)$-partitions, including the zero partition. Only $K\wt{v}\le N$ can contribute to degrees at most $N$.
\end{lemma}
\begin{proof}
Take any legal $d$-word. For a nonnull word, use its unique zero-separated parse
\[
 w=\iota(u_0)\,0\,\iota(u_1)\,0\cdots,
 \qquad u_0\supeq u_1\supeq\cdots,
\]
with absent trailing blocks padded by the null word; for the null word $w=\varnothing$, take all $u_i=\varnothing$ by convention. Fix the block $v=u_{K-1}$, allowing $v=\varnothing$ if there are fewer than $K$ nonnull blocks. The preceding $J=K-1$ word blocks $u_0,\ldots,u_{K-2}$ form a coupled $[J]$-family, each containing $v$. Its least admissible family is $v^{(J)}$, and its excess generating series is $C_{d-1}([J],v^{(J)};q)$. The succeeding blocks $u_K,u_{K+1},\ldots$ form a decreasing word chain bounded above by $v$, counted by $F(v;q)$. These choices are independent once $v$ is fixed: juxtaposing the shifted blocks with their zero separators uniquely reassembles the original legal word, and every such word has exactly one split at index $K-1$. By the recursive word-weight identity, the $J$ compulsory preceding copies and the fixed block contribute $K\wt{v}$; the remaining preceding and succeeding weights are precisely the two excess terms in the product. When $v=\varnothing$, the succeeding chain is necessarily null and the same argument applies. Every contribution therefore has weight at least $K\wt{v}$, proving the identity and the cutoff.
\end{proof}

\subsection{Conditioning a transverse word family}
Fix $k\ge2$, a grid $\Gamma$, a coupled $\Gamma$-family $\phi$ of legal $k$-words, and an integer $S\ge2$. Re-parse the members of $\phi$ at their outer zero separators; let $\phi_{\cdot,r}$ be the $\Gamma$-family of their $r$th lower-dimensional word blocks, padding absent blocks by $\varnothing$. By Proposition~\ref{prop:transverse} this is an eventually-null decreasing chain of coupled words. Form the extended grid $\Gamma^+=\Gamma\times[S-1]$. Fix the transverse word family $c$ at index $S-1$ (the $S$th member), consisting of legal $(k-1)$-words with $c\supeq\phi_{\cdot,S-1}$. Using only recursive word joins and weights, define
\begin{align}
 \psi(c)_{(\gamma,r)}&=c_\gamma\vee\phi_{\gamma,r}, \qquad 0\le r\le S-2,\label{eq:psi}\\
 \eps(c)&=\wt{c}-\wt{\phi_{\cdot,S-1}}
 +\sum_{r=0}^{S-2}\left(\wt{\psi(c)_{\cdot,r}}-\wt{\phi_{\cdot,r}}\right),\label{eq:eps}\\
 \Phi_r&=\phi_{\cdot,r}\qquad(r\ge S).\nonumber
\end{align}
By Proposition~\ref{prop:wordops}, the join of coupled nested word families is again coupled and nested; hence $\psi(c)$ is an admissible $\Gamma^+$-floor. Every summand of $\eps(c)$ is nonnegative. Since each of the $S-1$ upper transverse blocks in $\psi(c)$ contains $c$, the exact compulsory weight also gives
\begin{equation}\label{eq:transverse-forced}
 S\wt{c}\le\wt{\phi}+\eps(c).
\end{equation}
Consequently, when the allowed total weight is at most $N$, every affordable transverse conditioning state satisfies $S\wt{c}\le N$. This is the same forced-smallness mechanism as at the first word cut.

\begin{lemma}[Effective floor]\label{lem:effective}
For a coupled $\Gamma^+$-family $w$ of legal $(k-1)$-words, the two conditions
\[
 w_{\gamma,S-2}\supeq c_\gamma\quad\text{and}\quad
 w_{\gamma,r}\supeq\phi_{\gamma,r}\ (0\le r\le S-2)
\]
are equivalent to $w\supeq\psi(c)$.
\end{lemma}
\begin{proof}
If the first two conditions hold, nesting in $r$ gives $w_{\gamma,r}\supeq c_\gamma$ for every $r\le S-2$, while the floor condition gives $w_{\gamma,r}\supeq\phi_{\gamma,r}$, hence $w\supeq\psi(c)$. Conversely, $w\supeq\psi(c)$ immediately gives all floors, and at $r=S-2$ gives $w_{\gamma,S-2}\supeq c_\gamma$.
\end{proof}

\begin{proposition}[Recursive floor-conditioned decomposition]\label{prop:intermediate}
With the notation above,
\begin{equation}\label{eq:intermediate}
 C_k(\Gamma,\phi;q)
 =\sum_{c\supeq\phi_{\cdot,S-1}}
 q^{\eps(c)}C_{k-1}(\Gamma^+,\psi(c);q)\,G(c,S;\Phi;q).
\end{equation}
A column contributes through degree $E$ only if $\eps(c)\le E$, in which case
$\wt{c}\le E+\wt{\phi_{\cdot,S-1}}$.
\end{proposition}
\begin{proof}
An admissible family $u\supeq\phi$ is a coupled family of recursively parsed $k$-words. Proposition~\ref{prop:transverse} re-reads it uniquely as a chain of $(k-1)$-word families $u_{\cdot,r}$, nested in both the family index $\gamma$ and the block index $r$. Fix $c=u_{\cdot,S-1}$ and divide the remaining transverse words into the first $S-1$ and those at indices $r\ge S$. The upper words form a $\Gamma^+$-family with floor $\psi(c)$ by Lemma~\ref{lem:effective}. The later words form a chain bounded above by $c$, with floors $\Phi_r$, counted by $G(c,S;\Phi;q)$ after normalisation. Conversely, legal upper and lower word completions reassemble uniquely into $u$ by Proposition~\ref{prop:transverse}. The excess weight satisfies the exact equality
\[
 \wt{u}-\wt{\phi}
 =\bigl(\wt{\text{upper}}-\wt{\psi(c)}\bigr)
 +\eps(c)+\bigl(\wt{\text{lower}}-\sum_{r\ge S}\wt{\Phi_r}\bigr).
\]
This proves the identity. The size cutoff follows from $\eps(c)\ge\wt{c}-\wt{\phi_{\cdot,S-1}}$.
\end{proof}

\subsection{Normalised lower-chain recursion}
\begin{lemma}[Lower-chain transfer]\label{lem:lower}
For $X\supeq\Phi_r$,
\begin{equation}\label{eq:lower}
 G(X,r;\Phi;q)=\sum_{\Phi_r\subeq Y\subeq X}
 q^{\wt{Y}-\wt{\Phi_r}}G(Y,r+1;\Phi;q).
\end{equation}
After the floors vanish, $G=F$, $F(0;q)=1$, and for $X\ne0$,
\begin{equation}\label{eq:stationary}
 (1-q^{\wt{X}})F(X;q)=\sum_{Y\subsetneq X}q^{\wt{Y}}F(Y;q).
\end{equation}
These recurrences determine every truncated chain series uniquely.
\end{lemma}
\begin{proof}
Split each chain by its first member $Y_r=Y$. Normalising by $\sum_{j\ge r}\wt{\Phi_j}$ separates the first excess $\wt{Y}-\wt{\Phi_r}$ from the remaining normalised chain. Every exponent is nonnegative. Once the floor becomes identically zero, the first-member term $Y=X$ is $q^{\wt{X}}F(X;q)$ and can be moved to the left. Induction on $\wt{X}$ then determines $F(X;q)$; before the last nonzero floor, backward induction on the slice index determines $G$. The formal division by $1-q^{\wt{X}}$ is well defined because $\wt{X}>0$.
\end{proof}

Algorithms R and D are specified using this \emph{first-member} recursion, and the complexity analysis uses it throughout.

\section{Algorithm R: recursive row conditioning}\label{sec:alg}
The identities of Section~\ref{sec:conditioning} generalise Proposition~\ref{prop:ordinary-interior} and are iterations of the single transverse word-conditioning rule of Proposition~\ref{prop:intermediate}, and already apply at $k=2$. At this final row cut the constituent words are ordinary partition words. To make this a complete algorithm we need only specify their base case: counting coupled binary words by successive part positions. This section gives the base, the word-native state listing procedure, a complete recursion, and the correctness proof. The diagonal method will be introduced as an optional replacement in Section~\ref{sec:diagonal}.

\subsection{The integer-family base case}
A coupled family of ordinary-partition words can be read directly as sequences of part lengths: split each binary word at its zeroes and add one to the lengths of its runs of ones. A family $u=(u_\gamma)_{\gamma\in\Gamma}$ of such words gives a chain $x_0\succeq x_1\succeq\cdots$, where $x_r=(u_\gamma(r))_{\gamma\in\Gamma}$ is the tuple of $r$th part lengths, padded by zeros at absent positions. The family coupling and each word's weak decrease are precisely the inequalities on these integer tuples. A floor word family $\phi$ similarly gives a chain $f_r$. We regard a coupled integer tuple as a family of zero-dimensional partitions, with weight $|x|=\sum_\gamma x_\gamma$; this is simply the final parsing level of the same word algorithm.

For a fixed previous integer family $X$ and a floor tail $f_r,f_{r+1},\ldots$, define $G_0(X,r;f;q)$ to count remaining chains with exact excess over the floors. Its first-element recurrence is
\begin{equation}\label{eq:int-base}
G_0(X,r;f;q)=\sum_{f_r\preceq Y\preceq X}q^{|Y|-|f_r|}G_0(Y,r+1;f;q).
\end{equation}
After the floor vanishes, this reduces to the stationary chain series $F_0(X;q)$ with $F_0(0;q)=1$ and
\begin{equation}\label{eq:int-stationary}
(1-q^{|X|})F_0(X;q)=\sum_{0\preceq Y\prec X}q^{|Y|}F_0(Y;q),\qquad X\ne0.
\end{equation}
The first integer family is not bounded by a preceding member. Thus, for $E\ge0$,
\begin{equation}\label{eq:C1-base}
C_1(\Gamma,\phi;q)\big|_{\le E}=
\left[\sum_{\substack{x\succeq f_0\\ |x|-|f_0|\le E}}
q^{|x|-|f_0|}G_0(x,1;f;q)\right]_{\le E}.
\end{equation}
The total budget makes the first sum finite. The recurrence \eqref{eq:int-base} terminates because the floor eventually vanishes, and \eqref{eq:int-stationary} determines each remaining stationary state by induction on $|X|$.

For one ordinary partition, this reduces to choosing successive parts subject to weak decrease and a floor. For example, let
\[
H(b,m;L)=\#\{a\succeq b:|a|=m,\ a_0\le L\}.
\]
Writing $b'=(b_1,b_2,\ldots)$ gives, for $m>0$,
\begin{equation}\label{eq:Hrec}
H(b,m;L)=\sum_{s=\max(1,b_0)}^{\min(L,m)}H(b',m-s;s),
\end{equation}
with the appropriate empty-tail boundary case. This is the elementary weight-indexed recurrence underlying the hand count in Section~\ref{sec:example}; the coupled integer-family form preserves its compatibility constraints.

\subsection{Canonical word states and interval generation}
A legal partition word records its coordinate blocks uniquely. In the row algorithm a state is represented by its coupled word blocks, a compact floor word (or the unread suffix of that floor), a truncation budget, and the relevant position or conditioning depth. Containment and componentwise joins are computed blockwise. Word length plus one gives the weight of each nonempty encoded partition; the null partition has weight zero. A fixed part or row is represented by its own recursive word.

For exact enumeration of possible conditioning columns, centre states and chain successors, it suffices to generate monotone integer families in a constrained interval. The following elementary lemma makes this enumeration duplicate-free through word-native interval choices. Its formulation applies equally to finite product-grid ideal representations.
\begin{lemma}[Duplicate-free interval-word generation]\label{lem:wordinterval}
Let $P$ be a product of a fixed number $t$ of finite chains, equipped with its coordinatewise order. Let $f\preceq u$ be two monotone integer families on $P$, that is, weakly decreasing along the order of $P$; some upper bounds may instead be imposed by a finite nonnegative excess budget $B$. Suppose $\kappa$ is a nonnegative cost, zero at $f$, increasing by at least one whenever any admissible entry is increased by one, and suppose that its value and the cost changes required during interval generation can be evaluated in time polynomial in the encoded state, floor and budget sizes. Then all monotone families $x$ such that
\[
f\preceq x\preceq u,\qquad \kappa(x)\le B,
\]
are generated exactly once, without dead ends, by choosing entries along a fixed linear extension of $P$, each position after all positions below it, and enforcing the floor, ceiling and previously chosen predecessor bounds. For fixed $t$, the number of partial choices and elementary tests per output is polynomial in the finite budget and the encoded floor size. The statement also holds for a labelled coupled family whose indices form a product grid.
\end{lemma}
\begin{proof}
Fix a finite box containing every relevant position; such a box exists because each new nonzero entry costs at least one unit and the total excess is at most $B$. Process positions in an order in which every dominating predecessor appears before the dependent position. When selecting $x_z$, impose
\[
f_z\le x_z\le\min\{u_z,\,x_y:y\text{ is an already selected dominating predecessor of }z\}.
\]
If an upper ceiling is absent, the remaining budget gives a finite upper bound $f_z+B$. The interval is nonempty: the floor is monotone, every previously chosen predecessor satisfies $x_y\ge f_y\ge f_z$, and $u_z\ge f_z$ wherever specified. The cost of a partial choice is the cost of its floor completion, obtained by filling all as yet unchosen positions with their floor entries, and a partial choice is retained precisely when this cost is at most $B$. The floor completion respects all inequalities across the chosen--unchosen boundary, so it is an admissible completion of every retained partial choice. Thus no retained node is a dead end. Candidate values of $x_z$ are tried in increasing order from $f_z$; since each unit increase raises the cost by at least one, the first value whose floor completion exceeds $B$ ends the scan, so each retained node causes at most one rejected test. Every finished family has a unique sequence of chosen entries, giving duplicate-freeness. Each output corresponds to a path of length at most the number of positions in the finite box. Every retained partial node extends to an output and, for a fixed depth, is a prefix of at least one output; the number of partial nodes is therefore at most the number of positions times the number of outputs. All comparisons and word-address operations are polynomial in the represented size; cost updates have this property by the explicit evaluation hypothesis. The number of grid directions $t$ is fixed in the applications here, and the box has polynomially many relevant positions in $N$.
\end{proof}
\begin{remark}
The lemma is a blockwise interval procedure, \emph{not} the assertion that every word containing a floor word begins with that floor as a literal prefix. The latter is false. In the ordinary-partition case the entry choices take the familiar form $f_j\le x_j\le\min(u_j,x_{j-1})$. For coupled words, the ceiling includes the meet of already selected dominating members. The least remaining fill is always the floor, as in the proof. The monotone-cost condition holds for the compulsory excess and for the diagonal successor costs used later.
\end{remark}

\subsection{A complete recursive procedure}
All polynomials are truncated through degree at most the stated budget. By convention, any call to a truncated counting series with a negative budget returns zero; equivalently, successor sums may omit all choices whose cost exceeds the available budget. Memo tables include the constrained state, its unread floor tail when applicable, and the budget. The loops and successor sums in Algorithms~\ref{alg:row} and~\ref{alg:diagonal} are enumerated by Lemma~\ref{lem:wordinterval}. We write $\mathsf{Ch}(X;b)$ for the truncated stationary series $F(X;q)$ and $\mathsf{Low}(X,r;b,\Phi)$ for the floored series $G(X,r;\Phi;q)$; these work for any dimension of the members, including integers.

\Needspace{8\baselineskip}
\par\smallskip
\begin{center}
\captionsetup{hypcap=false}
\captionof{algorithm}{Algorithm R: row conditioning through integer families}\label{alg:row}
\end{center}
\begingroup\footnotesize
\begin{algorithmic}[1]
\Procedure{TopR}{$d,N;K,S_{d-1},\ldots,S_2$}
 \If{$N\le1$} \State \Return $(p_d(0),\ldots,p_d(N))=(1,\ldots,1)$ \EndIf
 \State $P\gets0$
 \ForAll{$v\text{ a $(d-1)$-partition with }K|v|\le N$}
  \State $E\gets N-K|v|$; $\Gamma\gets[K-1]$; $\phi\gets(v,\ldots,v)$
  \State $P\gets P+q^{K|v|}[\mathsf{FamR}(d-1,\Gamma,\phi,E)\,\mathsf{Ch}(v;E)]_{\le E}$
 \EndFor
 \State \Return coefficients of $P$ through degree $N$
\EndProcedure
\Procedure{FamR}{$k,\Gamma,\phi,E$}
 \If{$k=1$} \State \Return $\mathsf{Base1}(\Gamma,\phi,E)$ \EndIf
 \State $S\gets S_k$; $\Phi_r\gets\phi_{\cdot,r}$ for $r\ge S$; $H\gets0$
 \ForAll{$c\succeq\phi_{\cdot,S-1}$ with $\varepsilon(c)\le E$}
   \State $b\gets E-\varepsilon(c)$; $\psi\gets\psi(c)$
   \State $H\gets H+q^{\varepsilon(c)}[\mathsf{FamR}(k-1,\Gamma\times[S-1],\psi,b)\mathsf{Low}(c,S;b,\Phi)]_{\le b}$
 \EndFor
 \State \Return $[H]_{\le E}$
\EndProcedure
\Procedure{Base1}{$\Gamma,\phi,E$}
 \State Use \eqref{eq:C1-base} with $\mathsf{Low}$ given by \eqref{eq:int-base}
 \State \Return coefficients through degree $E$
\EndProcedure
\Procedure{Low}{$X,r;b,\Phi$}
 \If{every future floor vanishes} \State \Return $\mathsf{Ch}(X;b)$ \EndIf
 \State \Return $[\sum_{\Phi_r\preceq Y\preceq X}q^{|Y|-|\Phi_r|}\mathsf{Low}(Y,r+1;b-|Y|+|\Phi_r|,\Phi)]_{\le b}$
\EndProcedure
\Procedure{Ch}{$X;b$}
 \If{$X=0$} \State \Return $1$ \EndIf
 \State \Return $[\sum_{Y\prec X}q^{|Y|}\mathsf{Ch}(Y;b-|Y|)/(1-q^{|X|})]_{\le b}$
\EndProcedure
\end{algorithmic}
\endgroup
\par\smallskip

\begin{theorem}[Exactness and termination of Algorithm R]\label{thm:correct-row}
For every $d\ge3$, $N\ge0$ and admissible integer choices $K,S_{d-1},\ldots,S_2\ge2$, Algorithm~\ref{alg:row} terminates and returns exactly $p_d(0),\ldots,p_d(N)$.
\end{theorem}
\begin{proof}
The first-member chain recurrences \eqref{eq:lower}, \eqref{eq:int-base} and their stationary forms are valid by splitting a chain at its first member. After the floors vanish, the stationary recursion decreases the weight of the bounded member strictly; before that point the floor index advances until its finite support is exhausted. The base \eqref{eq:C1-base} therefore terminates and returns the exact coupled ordinary-partition completion series.

Proceed by induction on $k$. At $k=1$, the assertion follows from the preceding base case. At $k\ge2$, Proposition~\ref{prop:intermediate} decomposes every floored family uniquely into a conditioned column, a preceding $(k-1)$-family and a succeeding chain. By induction, the first is counted by $\mathsf{FamR}(k-1)$; the chain recurrence counts the second exactly. Their compulsory and residual weights add to the original excess, so truncation at $E-\varepsilon(c)$ discards only terms above budget. Lemma~\ref{lem:wordinterval} enumerates each admissible column and successor once; its strictly positive incremental costs guarantee finite search. The top identity in Lemma~\ref{lem:top} reconstructs exactly the unrestricted $d$-partitions. This completes the induction.
\end{proof}

The complexity proof starts from the first-member transfers specified here. Selective cache eviction preserves exactness, but arbitrary repeated evaluations after eviction are not covered by the bound.

\section{Algorithm D: the diagonal refinement}\label{sec:diagonal}
\subsection{The diagonal word reading: an interlacing block grammar}
The row-conditioning stages so far repeatedly re-read lower-dimensional blocks using recursive containment. At the diagonal bottom, the states are still partition words, but their successive blocks obey a different local rule: \emph{interlacing}. Let $\Gamma$ be a product grid with $M=\wt{\Gamma}$ members, and let $\phi$ be a coupled family of plane-partition words. To define their alternative diagonal reading, let $\pi(a,b)$ denote the part value decoded from one of these words. For $\ell\ge0$ put
\[
 d_\ell(\pi)(i)=\pi(i,i+\ell),\qquad
 d_{-\ell}(\pi)(i)=\pi(i+\ell,i).
\]
Each $d_\ell(\pi)$ is an ordinary partition, hence has its own word $W_1$ as in Definition~\ref{def:words}. The \emph{diagonal-block word} of $\pi$ is the eventually null, two-sided sequence
\begin{equation}\label{eq:diagonalword}
 \mathcal D(\pi)=\bigl(\ldots,W_1(d_{-2}(\pi)),W_1(d_{-1}(\pi))\,\big|\,
 W_1(d_0(\pi))\,\big|\,W_1(d_1(\pi)),W_1(d_2(\pi)),\ldots\bigr).
\end{equation}
For this two-sided notation we extend $W_1$ by $W_1(\varnothing)=\varnothing$. The null block must not be confused with the empty word $\epsilon$, which encodes a one-cell partition. The separators designate the central block and are not letters of $W_2$. A \emph{$\Gamma$-family of diagonal words} means $M$ such two-sided sequences, compared coordinatewise across the index grid. We write $\alpha=(\alpha_\gamma)_{\gamma\in\Gamma}$ for a family of diagonal partitions, whether or not the $W_1$ notation is displayed. Its total weight is the sum of the sizes of its blocks:
\[
 \wt\pi=\sum_{\ell\in\mathbb Z}\wt{d_\ell(\pi)}.
\]
This is a second reading of the same plane-partition word, \emph{not} its outer row-block reading $W_2(\pi)$: consecutive diagonal blocks are constrained by interlacing rather than containment. Each diagonal block is still an ordinary binary partition word, and all successor intervals, minimum completions and weight increments can be evaluated from its parsed part lengths and inherited floor words. The advantage of the diagonal reading is the stronger compulsory completion cost imposed on every diagonal block. The following classical geometric fact establishes the new grammar.

\begin{lemma}[Diagonal interlacing]\label{lem:interlace}
An eventually zero two-sided sequence $\delta_\ell$ of ordinary partitions is the diagonal sequence of a plane partition if and only if, for all $i\ge0$ and $\ell\ge0$,
\begin{equation}\label{eq:interlace}
 \delta_\ell(i)\ge\delta_{\ell+1}(i)\ge\delta_\ell(i+1),
 \qquad
 \delta_{-\ell}(i)\ge\delta_{-\ell-1}(i)\ge\delta_{-\ell}(i+1).
\end{equation}
A $\Gamma$-family lies above a floor $\phi$ precisely when these inequalities hold for each member, the diagonals are monotone in $\Gamma$, and $d_\ell(u_\gamma)\supeq d_\ell(\phi_\gamma)$ for all $\ell$. Given the centre $\lambda=d_0(u)$, the two sides are independent.
\end{lemma}
\begin{proof}
For $b\ge a$, the two inequalities comparing $\pi(a,b)$ with $\pi(a,b+1)$ and $\pi(a+1,b)$ translate into the first interlacing pair; the region $a\ge b$ is its transpose. These account for every adjacent inequality of a weakly decreasing array. Family ordering and floor constraints are coordinatewise and therefore hold precisely when they hold on every diagonal. No adjacent-cell inequality directly relates a strictly positive to a strictly negative diagonal once the centre is specified.
\end{proof}

In word language, Lemma~\ref{lem:interlace} gives three rules: consecutive nonnull blocks on either side interlace, hence decrease coordinatewise as one moves away from the centre; blocks from different family members remain ordered along $\Gamma$; and transposition of a plane partition exchanges the two sides. These conditions are both necessary and sufficient. Once the central block is prescribed, the two tails have no further mutual constraints, so their generating series multiply. This is the exact analogue of the separating-block operation at the upper recursive levels, with a different local grammar.

\begin{figure}[htbp]
\centering
\begin{tikzpicture}[font=\small,
  cell/.style={draw,minimum width=.7cm,minimum height=.7cm,inner sep=0pt,outer sep=0pt}]
\foreach \i/\j/\v in {0/0/6,0/1/5,0/2/3,1/0/4,1/1/3,1/2/2,2/0/2,2/1/1,2/2/1}{
  \pgfmathtruncatemacro{\diagoff}{\j-\i}
  \ifnum\diagoff=0 \def\shade{blue!25}\else\ifnum\diagoff=1 \def\shade{green!22}\else\ifnum\diagoff>1 \def\shade{gray!30}\else\def\shade{orange!22}\fi\fi\fi
  \node[cell,fill=\shade] at (.7*\j,-.7*\i) {$\v$};
}
\node[anchor=west] at (2.3,0) {centre: $(6,3,1)$};
\node[anchor=west] at (2.3,-.55) {positive side: $(5,2),(3)$};
\node[anchor=west] at (2.3,-1.10) {negative side: $(4,1),(2)$};
\node[anchor=west] at (2.3,-1.65) {$27=20+5+2$ (least completion + excesses)};
\end{tikzpicture}
\caption{An asymmetric plane partition with three distinct central diagonal entries. The blue cells form the central ordinary-partition word, the green and grey cells form the positive diagonals, and the orange cells the negative diagonals. The two sides interlace separately with the centre, but are not identical.}
\label{fig:diagonal-interlacing}
\end{figure}
In Figure~\ref{fig:diagonal-interlacing}, the central diagonal $(6,3,1)$ has least floor-free completion
of weight $6+3\cdot3+5\cdot1=20$. In the displayed partition,
the positive diagonals contribute five additional cells over that
completion, and the negative diagonals two, giving total weight
$20+5+2=27$. The displayed positive and negative tails are different choices of continuation from the same centre. Their generating series multiply because the choices are conditionally independent. For an empty floor the two series coincide by transposition; an asymmetric floor can make them different.

\subsection{The central word and its least geometric completion}
\begin{lemma}[Least completion of a centre]\label{lem:center}
Let $\lambda$ be a $\Gamma$-family of partitions with $\lambda_\gamma\supeq d_0(\phi_\gamma)$. The least $\Gamma$-family $u\supeq\phi$ of plane partitions with centre $d_0(u)=\lambda$ is
\begin{equation}\label{eq:center-completion}
 u_\gamma(a,b)=\max\{\lambda_\gamma(\max(a,b)),\phi_\gamma(a,b)\}.
\end{equation}
Its excess over the floor is
\begin{align}
 \Delta_0(\lambda)&=\sum_{\gamma\in\Gamma}\sum_{c\ge0}
 \bigl(\Psi_c^\gamma(\lambda_\gamma(c))-
 \Psi_c^\gamma(\phi_\gamma(c,c))\bigr),\label{eq:centercost}\\
 \Psi_c^\gamma(x)&=(2c+1)x+
 \sum_{a<c}\bigl(\phi_\gamma(a,c)-x\bigr)_+
 +\sum_{a<c}\bigl(\phi_\gamma(c,a)-x\bigr)_+.
 \label{eq:centerpsi}
\end{align}
For all admissible $x$, $\Psi_c^\gamma(x+1)-\Psi_c^\gamma(x)\ge1$.
\end{lemma}
\begin{proof}
Every admissible $u$ obeys $u(a,b)\ge u(\max(a,b),\max(a,b))=\lambda(\max(a,b))$ and $u\supeq\phi$. Their maximum is itself weakly decreasing in every coordinate and has the specified centre, proving minimality and admissibility, including the family order. Group cells by the largest of their two indices, $c=\max(a,b)$; this gives the expression for $\Psi$. At the floor centre $\lambda(c)=\phi(c,c)$ the right side yields the original floor weight. Finally,
\[
 \Psi_c^\gamma(x+1)-\Psi_c^\gamma(x)
 =(2c+1)-\#\{a<c:\phi_\gamma(a,c)>x\}
 -\#\{a<c:\phi_\gamma(c,a)>x\}\ge1.
\]
The estimate is uniform in the floor and in $\Gamma$.
\end{proof}

For an empty floor, the central block carries the position-dependent minimum cost
\begin{equation}\label{eq:centre-word-cost}
 \Delta_0(\lambda)=\sum_{\gamma\in\Gamma}
 \sum_{c\ge0}(2c+1)\lambda_\gamma(c).
\end{equation}
Thus one unit in part position $c$ of a central diagonal word forces $2c+1$ cells in its least completed plane partition. With a floor, the same statement is corrected \emph{exactly} by the positive-part terms in \eqref{eq:centerpsi}; they cannot be discarded or encoded by a constant weight on the alphabet letters alone.

\subsection{Legal side-word extensions and their excess cost}
\paragraph{Continuation states.} A positive-side state $(\ell,\alpha)$ consists of $\ell\ge0$ and a $\Gamma$-family of ordinary partitions $\alpha\supeq d_\ell(\phi)$. Let $\Cont(\ell,\alpha)$ be the eventually zero sequences $\alpha^{(m)}$, $m>\ell$, of $\Gamma$-families satisfying the positive interlacing inequalities, all floor constraints, and $\alpha^{(\ell)}=\alpha$. Their weight is $\sum_{m>\ell}\wt{\alpha^{(m)}}$. Write $L(\ell,\alpha)$ for the weight of a least continuation, and define
\[
 T(\ell,\alpha;q)=\sum_{\beta\in\Cont(\ell,\alpha)}
 q^{\wt{\beta}-L(\ell,\alpha)}.
\]
For the negative side use the transpose floor $\phi^\mathsf T$ and the same positive-side formulas.

\begin{lemma}[Least continuation, successors and transfer]\label{lem:side}
For every state $(\ell,\alpha)$ the following hold.
\begin{enumerate}[label=(\alph*),leftmargin=2.2em]
\item The pointwise least admissible continuation is
\begin{equation}\label{eq:least-side}
 \alpha^{*(m)}_\gamma(i)=
 \max\{\alpha_\gamma(i+m-\ell),\phi_\gamma(i,i+m)\},\quad m>\ell,
\end{equation}
with total weight
\begin{equation}\label{eq:L-side}
 L(\ell,\alpha)=\sum_{\gamma\in\Gamma}\sum_{j\ge0}
 \left[j\alpha_\gamma(j)+\sum_{a<j}\bigl(\phi_\gamma(a,j+\ell)-\alpha_\gamma(j)\bigr)_+\right].
\end{equation}
\item The possible first successors $\alpha'=\alpha^{(\ell+1)}$ are precisely the $\Gamma$-families of ordinary partitions in the interval
\begin{equation}\label{eq:beta}
 \beta^\gamma_i\le\alpha'^\gamma_i\le\alpha^\gamma_i,
 \qquad \beta^\gamma_i=\max\{\alpha^\gamma_{i+1},\phi_\gamma(i,i+\ell+1)\}.
\end{equation}
In particular $\beta$ itself is the least successor.
\item For the local functions
\begin{equation}\label{eq:f-side}
 f_i^\gamma(x)=(i+1)x+\sum_{a<i}
 \bigl(\phi_\gamma(a,i+\ell+1)-x\bigr)_+,
\end{equation}
the excess of a first successor is
\begin{equation}\label{eq:delta-side}
 \Delta(\alpha')=\sum_{\gamma,i}
 \left[f_i^\gamma(\alpha'^\gamma_i)-f_i^\gamma(\beta^\gamma_i)\right]\ge0,
\end{equation}
and the exact transfer is
\begin{equation}\label{eq:T-side}
 T(\ell,\alpha;q)=\sum_{\beta\subeq\alpha'\subeq\alpha}
 q^{\Delta(\alpha')}T(\ell+1,\alpha';q),\qquad T(\ell,0;q)=1.
\end{equation}
Every increment $f_i^\gamma(x+1)-f_i^\gamma(x)\ge1$.
\end{enumerate}
\end{lemma}
\begin{proof}
For (a), repeatedly applying the lower interlacing inequality forces $\alpha^{(m)}(i)\ge\alpha(i+m-\ell)$, and the floor forces $\alpha^{(m)}(i)\ge\phi(i,i+m)$. Their maximum in \eqref{eq:least-side} is a partition, because both arguments decrease in $i$. It is eventually zero because both the present partition and the floor are finitely supported. It obeys the upper interlacing inequality because the two arguments separately decrease with $m$; the lower inequality holds because $\alpha(i+m+1-\ell)$ is the shifted present bound and $\phi(i,i+m+1)\ge\phi(i+1,i+m+1)$. At $m=\ell+1$ it lies between $\alpha(i+1)$ and $\alpha(i)$ because $\alpha$ dominates the current diagonal of $\phi$. Monotonicity in $\Gamma$ is preserved by the maximum. Summing its entries and setting $j=i+m-\ell$ gives \eqref{eq:L-side}.

For (b), the inequalities at the first step are necessary. Conversely, any admissible family $\alpha'$ in the indicated coordinate interval is automatically a family of partitions: $\alpha'(i+1)\le\alpha(i+1)\le\beta_i\le\alpha'(i)$. It meets the next floor by definition of $\beta$, and (a) supplies at least one continuation. The family $\beta$ itself is admissible by the same argument.

For (c), first observe the compatibility of consecutive least continuations. The least continuation of $\beta$ at level $\ell+1$ equals that of $\alpha$ from level $\ell+2$ onward: by \eqref{eq:beta}, its value at $(i,m)$ is the maximum of $\alpha(i+m-\ell)$, $\phi(i+m-\ell-1,i+m)$ and $\phi(i,i+m)$, and the second floor term is bounded by the third along its column. Consequently
\[
 L(\ell,\alpha)=\wt{\beta}+L(\ell+1,\beta).
\]
A continuation whose first member is $\alpha'$ contributes an excess
\[
 \wt{\alpha'}+L(\ell+1,\alpha')
 -\wt{\beta}-L(\ell+1,\beta)
\]
plus the normalised excess of its tail. Applying \eqref{eq:L-side} to the first four terms yields exactly \eqref{eq:delta-side}. Moreover,
\[
 f_i^\gamma(x+1)-f_i^\gamma(x)
 =(i+1)-\#\{a<i:\phi_\gamma(a,i+\ell+1)>x\}\ge1.
\]
Thus $\Delta$ is nonnegative and locally additive, and partitioning the continuations by their first member proves \eqref{eq:T-side}.
\end{proof}

\begin{corollary}[Gap-word form of a floor-free side step]\label{cor:gapword}
For one plane partition and zero floor, let the current diagonal be $\alpha$, and put
\[
 g_i=\alpha_i-\alpha_{i+1}\quad(i\ge0).
\]
Every legal next diagonal $\alpha'$ is uniquely represented by finite nonnegative integers $x_i$ with
\begin{equation}\label{eq:gap-step}
 0\le x_i\le g_i,\qquad
 \alpha'_i=\alpha_{i+1}+x_i.
\end{equation}
Its successive gaps and its excess over the least successor are
\begin{equation}\label{eq:gap-update}
 g'_i=g_{i+1}+x_i-x_{i+1},\qquad
 \Delta(\alpha')=\sum_{i\ge0}(i+1)x_i.
\end{equation}
For a coupled family the same formulas apply memberwise, with the additional requirement that $\alpha'$ remain ordered over $\Gamma$.
\end{corollary}
\begin{proof}
The floor-free bounds in \eqref{eq:beta} are $\alpha_{i+1}\le\alpha'_i\le\alpha_i$, which are exactly \eqref{eq:gap-step}. Subtracting neighbouring entries yields the formula for $g'_i$; inserting the floor-free local function $f_i(x)=(i+1)x$ into \eqref{eq:delta-side} yields the cost. The successor is automatically a partition by the interlacing inequalities, and the variables are eventually zero.
\end{proof}

This is a position-decorated word extension, not a new fixed-alphabet letter weight: $i+1$ depends on the part position. The cheapest next block deletes the first part of $\alpha$ (all $x_i=0$). With nonzero floors its least extension is instead $\beta_i=\max\{\alpha_{i+1},\phi(i,i+\ell+1)\}$, and the exact extra cost is the corrected function in \eqref{eq:f-side}--\eqref{eq:delta-side}. In particular, the floor-dependent recursion of Lemma~\ref{lem:side} remains essential; the gap notation does not replace it.

\subsection{Rejoining the two word sides}
\begin{proposition}[The bottom-level family series]\label{prop:bottom}
Let $T^+$ and $T^-$ denote the transfers for $\phi$ and for the memberwise transpose $\phi^\mathsf T$. Then
\begin{equation}\label{eq:bottom}
 C_2(\Gamma,\phi;q)=\sum_{\lambda\supeq d_0(\phi)}
 q^{\Delta_0(\lambda)}T^+(0,\lambda;q)T^-(0,\lambda;q).
\end{equation}
Every nonzero diagonal $\alpha=d_\ell(u)$ of a completed family $u$ with $\wt{u}\le N$ has $\wt{\alpha}\le N$ and $|\ell|\le N-1$.
\end{proposition}
\begin{proof}
Given the centre, Lemma~\ref{lem:interlace} makes the two sides independent. Lemma~\ref{lem:center} supplies the least completion and its central excess; Lemma~\ref{lem:side} subtracts the least side costs and counts every additional excess exactly. A nonzero diagonal at distance $|\ell|$ includes an entry $u(0,|\ell|)\ge1$ after possible transposition, forcing at least $|\ell|+1$ cells along that row. Thus $|\ell|+1\le\wt{u}\le N$. Distinct family members occupy distinct labelled cells, so the same bounds apply to the summed diagonal weights.
\end{proof}

\begin{lemma}[Global budget and feasible diagonal states]\label{lem:size}
Every recursive family call reached from the top satisfies $|\phi|+E\le N$, where $\phi$ is its floor and $E$ its excess budget. Moreover, each diagonal centre and each side state reached by Algorithm D admits a complete coupled plane-partition extension above its inherited floor of total weight at most $N$.
\end{lemma}
\begin{proof}
At the top, $\phi=v^{(K-1)}$ and $E=N-K|v|$, so $|\phi|+E=N-|v|\le N$. If a family call chooses $c$ at compulsory cost $\varepsilon(c)$, then the upper recursive floor satisfies $|\psi(c)|\le|\phi|+\varepsilon(c)$ by the exact compulsory-cost identity. Hence $|\psi(c)|+E-\varepsilon(c)\le N$. Inducting on the recursive levels proves the first assertion. At a diagonal bottom call, every retained centre is compatible with its floor and the least centre completion from Lemma~\ref{lem:center}. Every subsequent side state is selected from a nonempty feasible interval; Lemma~\ref{lem:side} supplies a least future continuation at zero additional excess, while the opposite side can be completed by its least continuation. The total weight of the resulting family is its floor weight plus the compulsory and previously selected excesses, at most $|\phi|+E\le N$. This proves the second assertion.
\end{proof}

\subsection{Algorithm D and its exactness}
Algorithm D retains the outer conditioning identity \eqref{eq:top} and the floor-family recursion \eqref{eq:intermediate} at levels $k\ge3$. At $k=2$, it replaces $\mathsf{FamR}(2)$ and its integer-family base by the two-sided diagonal procedure $\mathsf{BotD}$ described below. It counts the centre using Proposition~\ref{prop:bottom}, and each side using Lemma~\ref{lem:side}. We write $[H]_{\le E}$ for truncation through degree $E$.
\begin{algorithm}[!htbp]
\footnotesize
\caption{The diagonal replacement at level two}\label{alg:diagonal}
\begin{algorithmic}[1]
\Procedure{FamD}{$k,\Gamma,\phi,E$}
 \If{$k=2$} \State \Return $\mathsf{BotD}(\Gamma,\phi,E)$ \EndIf
 \State \Return the level-$k$ sum \eqref{eq:intermediate} with recursive calls $\mathsf{FamD}(k-1)$
\EndProcedure
\Procedure{BotD}{$\Gamma,\phi,E$}
 \State $H\gets0$; initialise positive- and negative-side memo tables
 \ForAll{coupled centre families $\lambda\succeq d_0(\phi)$ with $\Delta_0(\lambda)\le E$}
  \State $b\gets E-\Delta_0(\lambda)$
  \State $H\gets H+q^{\Delta_0(\lambda)}[\mathsf{Side}^+(0,\lambda;b,\phi)\mathsf{Side}^-(0,\lambda;b,\phi^{\mathsf T})]_{\le b}$
 \EndFor
 \State \Return $[H]_{\le E}$
\EndProcedure
\Procedure{Side}{$\ell,\alpha;b,\phi$}
 \If{$\alpha=0$} \State \Return $1$ \EndIf
 \State $\beta\gets$ least feasible successor as in \eqref{eq:beta}
 \State \Return $[\sum_{\beta\preceq\alpha'\preceq\alpha,\,\Delta(\alpha')\le b}q^{\Delta(\alpha')}\mathsf{Side}(\ell+1,\alpha';b-\Delta(\alpha'),\phi)]_{\le b}$
\EndProcedure
\end{algorithmic}
\end{algorithm}
The top-level procedure $\mathsf{TopD}$ is $\mathsf{TopR}$ with each $\mathsf{FamR}$ call replaced by $\mathsf{FamD}$; the outer chain factor remains unchanged.
\begin{theorem}[Exactness of the diagonal refinement]\label{thm:correct-diag}
Algorithm D terminates and computes $p_d(0),\ldots,p_d(N)$ exactly, for every fixed $d\ge3$ and every admissible set of its conditioning parameters.
\end{theorem}
\begin{proof}
For a floored coupled plane-partition family, Proposition~\ref{prop:bottom} expresses its generating series exactly as the sum over centre diagonal families of the product of the two normalised interlacing side series. Lemmas~\ref{lem:center} and \ref{lem:side} prove the least completions and local excesses, so each coefficient is counted with its correct weight and without duplication. A nonzero diagonal at distance $\ell$ forces at least $\ell+1$ cells in a containing plane partition, so no nonzero side state is reachable for $\ell\ge N$. Hence side recursion terminates. At higher levels, exactness follows by the same induction on \eqref{eq:intermediate} as in Theorem~\ref{thm:correct-row}; the top identity then returns the required sequence.
\end{proof}

\section{Uniform complexity bounds}\label{sec:complexity}
This section proves Theorems~\ref{thm:row-main} and \ref{thm:diag-main}. Its bounds quantify the freedom left after conditioning: the outer blocks must be enumerated, the smaller residual word families must be stored or visited, and their constrained continuations must be evaluated. The conditioning depths are algorithmic choices, and balancing these competing costs gives the stated exponents. Exact legality and coupling remain mandatory in the recurrences; for an upper bound on the number of states we may count a larger class obtained by forgetting some coupling restrictions. This enlargement must not be confused with a simplification of the exact algorithm. The primary recurrences use first-member lower-chain transfers and an interval-state enumerator with polynomial overhead per state. The dimension $d$ is fixed throughout; constants $C_d$ do not depend on the growing conditioning parameters. We count operations on exact nonnegative integers and address bit complexity at the end.

\subsection{Partition growth and the interior-column invariant}
The known fixed-dimensional estimates \cite{BPA,Yeliussizov,Oganesyan} imply that for every fixed $k\ge1$ there is $b_k>0$ with
\begin{equation}\label{eq:growth-bound}
p_k(n)\le\exp\left(b_k n^{k/(k+1)}\right),\qquad n\ge0,
\end{equation}
after enlarging $b_k$ for small $n$.
\begin{lemma}[A uniform tuple bound]\label{lem:tuple}
For $k\ge1$, $M\ge1$ and $W\ge0$, the number of ordered $M$-tuples of $k$-partitions with combined weight at most $W$ does not exceed
\begin{equation}\label{eq:tuple-bound}
(W+1)^M\exp\left(b_k M^{1/(k+1)}W^{k/(k+1)}\right).
\end{equation}
\end{lemma}
\begin{proof}
Choose the individual weights $w_1,\ldots,w_M$ with total at most $W$; there are at most $(W+1)^M$ choices. By \eqref{eq:growth-bound} and concavity,
$\sum_iw_i^{k/(k+1)}\le M^{1/(k+1)}W^{k/(k+1)}$. Multiplying the bounds for each weight vector proves the assertion.
\end{proof}

\begin{lemma}[Forced smallness of an interior column]\label{lem:strongsize}
Suppose a level-$k$ family call, $k\ge2$, has floor $\phi$ and excess budget $E$, with $|\phi|+E\le N$. If a conditioned column $c$ at depth $S=S_k$ satisfies $\varepsilon(c)\le E$, then
\begin{equation}\label{eq:small-column}
|c|\le N/S.
\end{equation}
The same bound applies to every lower-chain state reached below $c$.
\end{lemma}
\begin{proof}
Every one of the $S-1$ preceding effective-floor slices contains $c$, so $|\psi(c)|\ge(S-1)|c|$. The exact compulsory-weight identity is
\begin{equation}\label{eq:strongsize-identity}
|\phi|+\varepsilon(c)
=|\psi(c)|+|c|+\sum_{r\ge S}|\phi_{\cdot,r}|
\ge S|c|.
\end{equation}
The hypothesis gives $|\phi|+\varepsilon(c)\le|\phi|+E\le N$. Every lower-chain state is contained in $c$, so its weight is also bounded by $N/S$.
\end{proof}

\subsection{The row-only state count and its uniform cost}
Let $K\ge2$ denote the outer conditioning depth and $S_k\ge2$ the depth at each family level $2\le k\le d-1$. Set
\begin{equation}\label{eq:Gamma}
M_k=(K-1)\prod_{j=k+1}^{d-1}(S_j-1),\qquad1\le k\le d-1.
\end{equation}
Thus $M_{d-1}=K-1$, while $M_1=(K-1)\prod_{j=2}^{d-1}(S_j-1)$ is the number of coupled integer entries at the base. The grid of level-$k$ families has $M_k$ members and $d-1-k$ grid directions beyond the initial chain direction.

\begin{lemma}[Intermediate row-state count]\label{lem:row-count}
For each $2\le k\le d-1$ set
\begin{equation}\label{eq:A-k}
A_k=(N+1)^{M_k}\exp\left(b_{k-1}M_k^{1/k}(N/S_k)^{(k-1)/k}\right).
\end{equation}
Every conditioned column and every subsequent lower-chain family at level $k$ belongs to a population of at most $A_k$ candidates. The first-member-chain implementation uses at most $(N+1)^{C_d}A_k^2$ arithmetic operations beyond the recursive calls to level $k-1$, per level-$k$ family call, up to a polynomial factor whose exponent depends only on $d$.
\end{lemma}
\begin{proof}
Lemma~\ref{lem:strongsize} bounds the total weight of each column and lower-family state by $N/S_k$. Forget the coupling and floors. The remaining ordered $M_k$-tuple of $(k-1)$-partitions is bounded by Lemma~\ref{lem:tuple}, giving $A_k$. A lower-chain memo key consists of a possible family, a floor-position index and a truncation budget; only $O(N+1)$ positions before the last nonzero floor and at most $N+1$ budgets are relevant. Each state has at most $A_k$ candidate successors. The interval-word lemma gives polynomial overhead in $N$ per candidate. There are at most $N^{O_d(1)}$ available coordinate positions in the fixed-dimensional cell representation, and each state of weight at most $N$ admits encoding, joins, comparisons and truncated coefficient operations in polynomial time. The number of recursive upper-family calls is bounded by $A_k$.
\end{proof}

\begin{lemma}[The integer-family base cost]\label{lem:base-cost}
With $M=M_1$ coupled integer entries at the base, Algorithm R requires at most
\begin{equation}\label{eq:base-cost}
Q_1(N,M)\le(N+1)^{2M+C_d}
\end{equation}
arithmetic operations per base-family call.
\end{lemma}
\begin{proof}
Each integer entry is between 0 and $N$ in a state contributing at weight at most $N$. Ignoring order constraints, there are at most $(N+1)^M$ possible states. The first-member or stationary memo carries a position and truncation index, adding a polynomial factor. Each state visits at most $(N+1)^M$ successors. The state-listing and integer polynomial operations contribute only polynomial overhead. This proves \eqref{eq:base-cost} after enlarging $C_d$.
\end{proof}

\begin{proposition}[Uniform row-only operation bound]\label{prop:row-uniform}
For $2\le K,S_k\le N$ and fixed $d\ge3$, Algorithm R uses at most $\exp(\Lambda_d^{\mathrm R})$ arithmetic operations, where
\begin{align}
\Lambda_d^{\mathrm R}\le C_d\Bigg[&
 (N/K)^{(d-1)/d}
 +\sum_{k=2}^{d-1}M_k^{1/k}(N/S_k)^{(k-1)/k}\nonumber\\
 &+\left(K+\sum_{k=1}^{d-1}M_k\right)\log(N+2)
\Bigg].\label{eq:row-uniform}
\end{align}
\end{proposition}
\begin{proof}
By the tuple bound, the number of possible outer slices $v$ with $|v|\le N/K$ is at most $\,(N+1)\exp(b_{d-1}(N/K)^{(d-1)/d})$. The outer stationary chain memo visits no more shapes, and its successors lie within the same population. At each level $k\ge2$, Lemma~\ref{lem:row-count} gives a cost bounded by $A_k^2Q_{k-1}(N+1)^{C_d}$, where $Q_{k-1}$ bounds the next family call. Terminate the fixed number of levels with Lemma~\ref{lem:base-cost}. Taking logarithms, and absorbing polynomial factors, gives \eqref{eq:row-uniform}. The constants are uniform in the conditioning parameters because they never enter the exponent of the polynomial state-processing overhead.
\end{proof}

\begin{proof}[Proof of Theorem~\ref{thm:row-main}]
For $N\le1$ the algorithm returns the initial coefficients directly. For $N\ge2$, put $a=1/(d+1)$ and $\theta=(d-1)/(d+1)$, and choose
\[
K=S_2=\cdots=S_{d-1}=\max\{2,\lceil N^{a}\rceil\},
\]
so that $2\le K,S_k\le N$. Equation~\eqref{eq:Gamma} gives
$M_k=O_d(N^{(d-k)a})$. The outer term in
\eqref{eq:row-uniform} is $O_d(N^\theta)$, as is each
intermediate term, because
\[
\frac{d-k}{k(d+1)}+\frac{d(k-1)}{k(d+1)}
=\frac{d-1}{d+1}.
\]
Finally, $K+\sum_{k=1}^{d-1}M_k=O_d(N^\theta)$.
The last term in \eqref{eq:row-uniform} is therefore
$O_d(N^\theta\log(N+2))$, proving the stated bound.
Exactness is Theorem~\ref{thm:correct-row}.
\end{proof}

\paragraph{Solid partitions.}
For $d=3$, $M_2=K-1$ and $M_1=(K-1)(S_2-1)$.
The three leading terms in \eqref{eq:row-uniform} are
\[
(N/K)^{2/3},\qquad(KN/S_2)^{1/2},\qquad KS_2\log(N+2).
\]
Taking $K\asymp S_2\asymp N^{1/4}$ makes these respectively
$O(\sqrt N)$, $O(\sqrt N)$ and $O(\sqrt N\log(N+2))$.
The row method is exact at any admissible depths; these
balanced growing depths give the simplified uniform bound.

\subsection{The diagonal state's square completion weight}\label{subsec:square-weight}

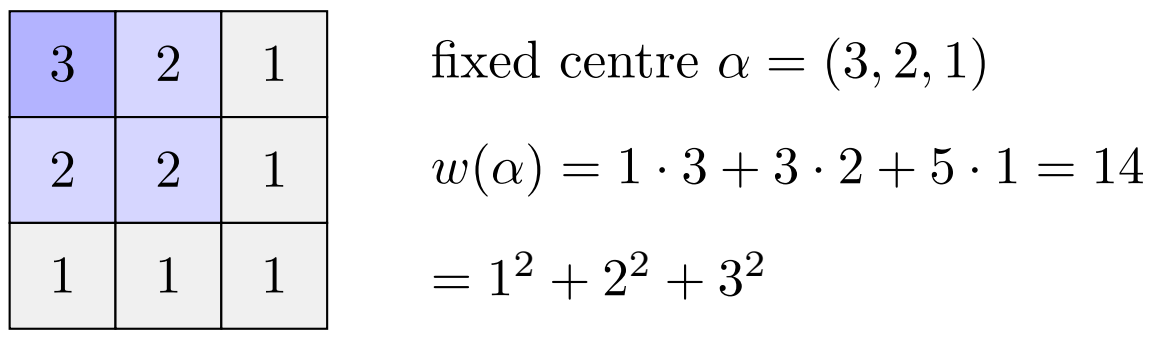
\begin{figure}[htbp]
\centering
\begin{tikzpicture}[font=\small,
  cell/.style={draw,minimum width=.7cm,minimum height=.7cm,inner sep=0pt,outer sep=0pt}]
\foreach \i in {0,1,2}{
  \foreach \j in {0,1,2}{
    \pgfmathtruncatemacro{\k}{max(\i,\j)}
    \pgfmathtruncatemacro{\val}{3-\k}
    \ifnum\k=0 \def\shade{blue!30}\else\ifnum\k=1 \def\shade{blue!16}\else\def\shade{gray!12}\fi\fi
    \node[cell,fill=\shade] at (.7*\j,-.7*\i) {$\val$};
  }
}
\node[anchor=west] at (2.3,0) {fixed centre $\alpha=(3,2,1)$};
\node[anchor=west] at (2.3,-.7) {$w(\alpha)=1\cdot3+3\cdot2+5\cdot1=14$};
\node[anchor=west] at (2.3,-1.4) {$=1^2+2^2+3^2$};
\end{tikzpicture}
\caption{The least plane-partition completion of the central diagonal word with decoded parts $(3,2,1)$. The rings $\max(a,b)=j$ contain $2j+1$ positions. Successive part differences convert this compulsory weight into a sum of square costs.}
\label{fig:diagonal-square-completion}
\end{figure}

For an ordinary partition $\alpha=(\alpha_0,\alpha_1,\ldots)$, define its \emph{completion weight} (Figure~\ref{fig:diagonal-square-completion})
\begin{equation}\label{eq:squareweight}
 w(\alpha)=\sum_{j\ge0}(2j+1)\alpha_j.
\end{equation}

\begin{lemma}[All diagonal states obey the square-weight constraint]\label{lem:squares}
Let $u$ be a plane partition, and let $\alpha=d_\ell(u)$ be any diagonal, central or noncentral. Then
\begin{equation}\label{eq:floorbound}
 \wt u\ge w(\alpha).
\end{equation}
If $\alpha=(\alpha_\gamma)_{\gamma\in\Gamma}$ is any centre or side state occurring in a bottom call whose total completion weight is at most $N$, then
\begin{equation}\label{eq:familyweight}
 \sum_{\gamma\in\Gamma}w(\alpha_\gamma)\le N.
\end{equation}
\end{lemma}
\begin{proof}
By transposition it suffices to treat $\ell\ge0$. For every pair $(a,b)$ put $j=\max(a,b)$. Coordinate monotonicity gives
\[
 u(a,b)\ge u(j,j)\ge u(j,j+\ell)=\alpha_j.
\]
There are exactly $2j+1$ pairs $(a,b)\in\NN^2$ with $\max(a,b)=j$, so summing proves \eqref{eq:floorbound}. Summation over members proves the family inequality for every complete $\Gamma$-family $u$. By Lemma~\ref{lem:size} and the least-continuation construction in Lemma~\ref{lem:side}, every memoised side state has an admissible complete extension of total weight at most $N$. Therefore \eqref{eq:familyweight} applies to \emph{every} reachable side state, whether its floor is zero or not, not merely to the centres.
\end{proof}

\begin{remark}[Distance-sensitive completion cost]\label{rem:distanceword}
The same geometric argument gives a sharper restriction for a diagonal $\alpha=d_\ell(u)$ at a specified distance $\ell\ge0$ from the centre:
\begin{equation}\label{eq:distance-weight}
 \wt u\ge w_\ell(\alpha)
 :=\sum_{j\ge0}(2j+\ell+1)\alpha_j
 =\sum_{r\ge1}r(r+\ell)m_r,
 \qquad m_r=\alpha_{r-1}-\alpha_r.
\end{equation}
Indeed, coordinate monotonicity forces $u(a,b)\ge\alpha_{\max\{a,b-\ell\}}$ for every $(a,b)$; the least array satisfying this condition has precisely $2j+\ell+1$ cells in the set where $\max\{a,b-\ell\}=j$. Transposition gives the same bound on the negative side. It is sharp without floors. Hence admissible side-word blocks farther from the centre satisfy more restrictive position-dependent costs; this may be used as an optional level-aware search bound. We use only $w_0=w$ in the uniform complexity proof, so no improvement of its leading exponent is claimed here.
\end{remark}

\begin{lemma}[Uniform count of square-weighted families]\label{lem:squarecount}
For the explicit constant $C=1+\frac{\sqrt\pi}{2}\zeta(3/2)$, the number of $M$-tuples of ordinary partitions satisfying $\sum_{\gamma=1}^M w(\alpha_\gamma)\le N$ is at most
\begin{equation}\label{eq:squarecount}
 D(M,N):=\exp\bigl(C M^{2/3}N^{1/3}\bigr)
\end{equation}
for all $M,N\ge1$. Coupling and floors can only reduce this number.
\end{lemma}
\begin{proof}
Let $m_r=\alpha_{r-1}-\alpha_r$ for $r\ge1$, with the convention that $\alpha_j=0$ eventually. Then
\begin{equation}\label{eq:squaresidentity}
 w(\alpha)=\sum_{r\ge1}m_r\sum_{j=0}^{r-1}(2j+1)
 =\sum_{r\ge1}r^2m_r.
\end{equation}
Thus ordinary partitions weighted by $w$ are exactly partitions into squares. Equivalently, if $g_i$ denotes the gap coordinates of Corollary~\ref{cor:gapword}, then $w(\alpha)=\sum_{i\ge0}(i+1)^2g_i$: the \emph{gaps} of a central diagonal word have square completion costs. The generating function for $M$ independent such partitions is
\[
 Z_M(q)=\prod_{r\ge1}(1-q^{r^2})^{-M}.
\]
Let $A_M(N)$ count tuples of total square weight at most $N$. Positivity gives, for every $t>0$,
\[
 A_M(N)\le e^{tN} Z_M(e^{-t}).
\]
Now
\begin{align*}
 \log Z_M(e^{-t})
 &=M\sum_{h\ge1}\frac1h\sum_{r\ge1}e^{-thr^2}\\
 &\le M\sum_{h\ge1}\frac1h\int_0^\infty e^{-thx^2}\,\dd x
 =\frac{\sqrt\pi}{2}\,\zeta(3/2)M t^{-1/2}.
\end{align*}
Choose $t=(M/N)^{2/3}>0$. The two terms sum to $\bigl(1+\frac{\sqrt\pi}{2}\zeta(3/2)\bigr) M^{2/3}N^{1/3}$, yielding the stated bound for all positive $M,N$, without an assumption that $M\le N$. This argument is uniform when $M$ grows with $N$.
\end{proof}

\subsection{The diagonal bottom cost and Theorem~\ref{thm:diag-main}}
\begin{lemma}[Uniform diagonal-bottom cost]\label{lem:cost-bottom}
For fixed $d$ and a bottom call with $M=M_2$ coupled plane-partition words, the diagonal procedure uses at most
\begin{equation}\label{eq:bottom-cost}
Q_2^{\mathrm D}(N,M)\le(N+1)^{C_d}\exp(C_dM^{2/3}N^{1/3})
\end{equation}
arithmetic operations.
\end{lemma}
\begin{proof}
Every listed centre and every reachable side state is a tuple of $M$ partitions satisfying \eqref{eq:familyweight}. By Lemma~\ref{lem:squarecount}, the number of such tuples is at most $D=\exp(CM^{2/3}N^{1/3})$. A side memo has at most $(N+1)^2D$ keys, including the diagonal index and remaining budget. Every state enumerates at most $D$ successors. Lemma~\ref{lem:wordinterval}, applied to the floor--ceiling intervals of Lemmas~\ref{lem:center} and \ref{lem:side}, gives polynomial overhead per state, as does convolution of the side coefficient vectors. The total is at most $(N+1)^{C_d}D^2$, absorbed by the stated bound.
\end{proof}

\begin{proposition}[Uniform diagonal operation bound]\label{prop:diag-uniform}
For $2\le K,S_3,\ldots,S_{d-1}\le N$, Algorithm D uses at most $\exp(\Lambda_d^{\mathrm D})$ arithmetic operations, where
\begin{align}
\Lambda_d^{\mathrm D}\le C_d\Bigg[&
 (N/K)^{(d-1)/d}
 +\sum_{k=3}^{d-1}M_k^{1/k}(N/S_k)^{(k-1)/k}\nonumber\\
 &+M_2^{2/3}N^{1/3}
 +\left(K+\sum_{k=2}^{d-1}M_k\right)\log(N+2)
\Bigg].\label{eq:diag-uniform}
\end{align}
\end{proposition}
\begin{proof}
The outer-slice population has logarithm at most $O_d((N/K)^{(d-1)/d}+\log(N+2))$, as before. At level $k\ge3$, Lemma~\ref{lem:row-count} applies unchanged, giving a factor $A_k^2(N+1)^{C_d}$ in the operation count. At the bottom, use Lemma~\ref{lem:cost-bottom} instead of the integer-family base. Multiplication through the fixed number of levels and taking logarithms gives \eqref{eq:diag-uniform}.
\end{proof}

\begin{proof}[Proof of Theorem~\ref{thm:diag-main}]
Correctness is Theorem~\ref{thm:correct-diag}. For $N\le1$ the algorithm returns the initial coefficients directly. For $N\ge2$, put $a=1/(d+1)$ and choose
\[
K=S_3=\cdots=S_{d-1}=\max\{2,\lceil N^a\rceil\},
\]
with only $K$ present when $d=3$. We have $M_k=O_d(N^{(d-k)a})$. In \eqref{eq:diag-uniform}, the outer term has exponent $(d-1)(1-a)/d=(d-1)/(d+1)$. Each intermediate term has exponent
\[
\frac{d-k}{k(d+1)}+\frac{d(k-1)}{k(d+1)}=\frac{d-1}{d+1}.
\]
The bottom term has exponent
\[
\frac{2(d-2)}{3(d+1)}+\frac13=\frac{d-1}{d+1}.
\]
Finally, $K+\sum_{k=2}^{d-1}M_k=O_d(N^{(d-2)/(d+1)})$, so the logarithmic overhead is lower order. The asserted bound follows.
\end{proof}

\begin{remark}[Asymptotic balancing versus finite-size choices]\label{rem:finite}
The parameters used to attain the asymptotic operation bounds are not asserted to minimise finite-size running times. Nearest-integer rounding of $N^{1/4}$ gives $K=3$ for approximately $40\le N<150$, but this does not account for implementation-dependent constants. The row and diagonal algorithms remain exact for every admissible parameter choice. With $K$ fixed independently of $N$, however, the outer term remains $O(N^{(d-1)/d})$: for $d=3$ and $K=3$ this argument gives only $\exp(O(N^{2/3}))$, rather than the $\exp(O(\sqrt N))$ bound attained using growing $K$. At $N=100$, $N^{1/4}\approx3.16$, whose nearest integer is $3$, whereas $\lceil N^{1/4}\rceil=4$.
\end{remark}

\paragraph{Solid partitions.}
For $d=3$, the two dominant logarithmic-cost contributions are $(N/K)^{2/3}$ and $K^{2/3}N^{1/3}$. Their balance is $K\asymp N^{1/4}$, and both become $O(\sqrt N)$. At finite $N$ this is a guide to parameter selection, not a statement that the integer nearest $N^{1/4}$ minimises actual measured runtime.

\paragraph{Bit complexity and scope.}
Every coefficient counts families of at most $N$ cells in a fixed-dimensional labelled product order. A generous upper bound of $\exp(O_d(N\log(N+2)))$ on intermediate counts implies that all exact integers manipulated have polynomially many bits in $N$ for fixed $d$. Standard exact-integer addition and multiplication increase the arithmetic-operation counts by at most a polynomial in $N$, absorbed into the displayed stretched exponentials. No estimate uniform in growing $d$ is claimed. The proof assumes that each constrained state is evaluated only a bounded number of times; an arbitrary policy of cache eviction and later recomputation is not thereby bounded.

\section{Discussion and limitations}\label{sec:discussion}
\paragraph{Which ingredients are new?} Nested slice descriptions of multidimensional partitions are elementary; diagonal interlacing is classical \cite{OkounkovReshetikhin}; and the asymptotic growth of square-partition coefficients belongs to the classical weighted-partition literature \cite{Wright}. The principal contribution is a method of exact enumeration by recursive interior conditioning of partition words, using transverse re-parsing, effective floors and nonnegative compulsory excess. Algorithm R realises the method through successive row conditioning; Algorithm D incorporates an interlacing diagonal transfer at the final level and achieves the stronger bound \eqref{eq:headline}. Both count constrained completions without enumerating finished partitions. Residual-floor state sharing (Appendix~\ref{app:residual}) is an additional implementation theorem.

\paragraph{Limitations and an open question.}
Both theorems are upper bounds for fixed $d$; no bound uniform in $d$ is claimed, and the balanced choice of parameters controls the derived exponent, not necessarily the runtime of every exact counting algorithm. Memory usage is not separately optimised in either theorem, and neither theorem resolves the asymptotic relationship between $p_d(n)$ and MacMahon's numbers. The logarithm in Theorem~\ref{thm:row-main} arises at the final row level, where Algorithm~R lists coupled integer families. The diagonal bottom avoids it because its reachable blocks obey square completion weights (Section~\ref{subsec:square-weight}). Whether a row-only algorithm can avoid the logarithm remains open. A persistent coefficient-state cache across separate target weights is also conceivable, but it requires retaining the relevant transfer graph and has not been measured here.

\section{Computed coefficients}\label{sec:computation}
\subsection{New coefficients}
Tables~\ref{tab:p3new}--\ref{tab:p7new} list the coefficients computed with Algorithm D beyond the ranges previously tabulated in the OEIS \cite{OEISSolid,OEIS4,OEIS5,OEIS6} and by Graffeo, Monavari, Moschetti and Ricolfi \cite{GMMR2026,GMMRData}: $p_3(n)$ for $73\le n\le100$, $p_4(n)$ for $41\le n\le50$, $p_5(n)$ for $33\le n\le45$, $p_6(n)$ for $33\le n\le40$ and $p_7(n)$ for $33\le n\le35$. The five-dimensional table also shows $p_5(31)$ and $p_5(32)$, which GMMR's table contains; see Section~\ref{sec:gmmr-discrepancy}. Spaces in large integers are grouping separators only. Throughout, $p_d(n)$ uses the array-dimension convention of Definition~\ref{def:kpartition}.

\begingroup
\small
\begin{longtable}{@{}r@{\qquad}r@{}}
\caption{Solid-partition coefficients $p_3(73),\ldots,p_3(100)$, beyond the range $n\le72$ of \cite{OEISSolid}.}\label{tab:p3new}\\
\toprule
$n$ & $p_3(n)$\\\midrule
\endfirsthead
\multicolumn{2}{@{}l}{\footnotesize Table~\thetable\ (continued)}\\
\toprule
$n$ & $p_3(n)$\\\midrule
\endhead
\midrule
\multicolumn{2}{r}{\footnotesize Continued on next page}\\
\endfoot
\bottomrule
\endlastfoot
73 & 5\,546\,716\,388\,521\,137\,692\\
74 & 8\,866\,970\,207\,488\,458\,458\\
75 & 14\,152\,754\,466\,500\,976\,843\\
76 & 22\,555\,039\,558\,627\,198\,898\\
77 & 35\,891\,640\,262\,751\,219\,089\\
78 & 57\,029\,563\,918\,064\,668\,233\\
79 & 90\,484\,360\,954\,223\,277\,574\\
80 & 143\,358\,415\,090\,239\,954\,282\\
\addlinespace[3pt]
81 & 226\,807\,963\,517\,255\,531\,305\\
82 & 358\,333\,889\,257\,123\,661\,244\\
83 & 565\,354\,399\,842\,777\,981\,841\\
84 & 890\,770\,070\,376\,764\,164\,789\\
85 & 1\,401\,621\,831\,357\,922\,117\,069\\
86 & 2\,202\,544\,821\,764\,780\,402\,353\\
87 & 3\,456\,650\,041\,671\,482\,513\,941\\
88 & 5\,417\,897\,027\,736\,772\,529\,969\\
89 & 8\,481\,223\,512\,538\,764\,767\,738\\
90 & 13\,260\,084\,406\,108\,380\,778\,701\\
\addlinespace[3pt]
91 & 20\,706\,253\,003\,113\,498\,210\,245\\
92 & 32\,294\,715\,698\,836\,959\,449\,895\\
93 & 50\,308\,718\,871\,164\,926\,060\,930\\
94 & 78\,278\,746\,849\,363\,887\,476\,119\\
95 & 121\,657\,842\,187\,627\,811\,973\,789\\
96 & 188\,859\,428\,795\,517\,348\,420\,742\\
97 & 292\,850\,580\,990\,579\,853\,624\,589\\
98 & 453\,595\,525\,513\,443\,249\,884\,458\\
99 & 701\,799\,328\,644\,219\,387\,872\,150\\
100 & 1\,084\,637\,900\,920\,436\,453\,116\,970\\
\end{longtable}
\endgroup

\begingroup
\small
\begin{longtable}{@{}rr@{}}
\caption{Four-dimensional coefficients $p_4(41),\ldots,p_4(50)$, beyond the range $n\le40$ of \cite{OEIS4}.}\label{tab:p4new}\\
\toprule
$n$ & $p_4(n)$\\\midrule
\endfirsthead
\multicolumn{2}{@{}l}{\footnotesize Table~\thetable\ (continued)}\\
\toprule
$n$ & $p_4(n)$\\\midrule
\endhead
\midrule
\multicolumn{2}{r}{\footnotesize Continued on next page}\\
\endfoot
\bottomrule
\endlastfoot
41 & 234\,310\,483\,718\,511\\
42 & 463\,193\,286\,245\,475\\
43 & 912\,703\,865\,994\,165\\
44 & 1\,792\,766\,221\,191\,810\\
45 & 3\,510\,520\,778\,098\,410\\
46 & 6\,853\,314\,878\,029\,174\\
47 & 13\,339\,437\,043\,652\,531\\
48 & 25\,888\,556\,874\,608\,490\\
49 & 50\,099\,922\,932\,879\,425\\
50 & 96\,682\,990\,637\,991\,290\\
\end{longtable}
\endgroup

\begingroup
\small
\begin{longtable}{@{}rr@{}}
\caption{Five-dimensional coefficients $p_5(31),\ldots,p_5(45)$. The OEIS lists $p_5$ through $n=30$ \cite{OEIS5}, and GMMR's table through $n=32$ \cite{GMMRData}; for $p_5(32)$ see Section~\ref{sec:gmmr-discrepancy}.}\label{tab:p5new}\\
\toprule
$n$ & $p_5(n)$\\\midrule
\endfirsthead
\multicolumn{2}{@{}l}{\footnotesize Table~\thetable\ (continued)}\\
\toprule
$n$ & $p_5(n)$\\\midrule
\endhead
\midrule
\multicolumn{2}{r}{\footnotesize Continued on next page}\\
\endfoot
\bottomrule
\endlastfoot
31 & 9\,385\,148\,192\,238\\
32 & 21\,684\,370\,682\,500\\
33 & 49\,883\,015\,219\,960\\
34 & 114\,274\,883\,300\,637\\
35 & 260\,750\,207\,205\,695\\
36 & 592\,720\,681\,792\,084\\
37 & 1\,342\,441\,905\,336\,952\\
38 & 3\,029\,860\,697\,217\,176\\
39 & 6\,815\,333\,428\,229\,786\\
40 & 15\,280\,540\,462\,255\,553\\
41 & 34\,152\,474\,835\,110\,987\\
42 & 76\,098\,745\,464\,055\,291\\
43 & 169\,059\,821\,812\,420\,048\\
44 & 374\,492\,239\,761\,142\,250\\
45 & 827\,205\,156\,290\,467\,825\\
\end{longtable}
\endgroup

\begingroup
\small
\begin{longtable}{@{}rr@{}}
\caption{Six-dimensional coefficients $p_6(33),\ldots,p_6(40)$, beyond the range $n\le32$ of GMMR's table \cite{GMMRData}; the OEIS lists $p_6$ through $n=21$ \cite{OEIS6}.}\label{tab:p6new}\\
\toprule
$n$ & $p_6(n)$\\\midrule
\endfirsthead
\multicolumn{2}{@{}l}{\footnotesize Table~\thetable\ (continued)}\\
\toprule
$n$ & $p_6(n)$\\\midrule
\endhead
\midrule
\multicolumn{2}{r}{\footnotesize Continued on next page}\\
\endfoot
\bottomrule
\endlastfoot
33 & 1\,623\,932\,959\,099\,444\\
34 & 4\,120\,190\,410\,516\,313\\
35 & 10\,407\,310\,523\,992\,865\\
36 & 26\,177\,828\,638\,944\,789\\
37 & 65\,583\,936\,794\,843\,761\\
38 & 163\,688\,518\,996\,765\,112\\
39 & 407\,077\,835\,092\,630\,492\\
40 & 1\,008\,905\,774\,197\,319\,306\\
\end{longtable}
\endgroup

\begingroup
\small
\begin{longtable}{@{}rr@{}}
\caption{Seven-dimensional coefficients $p_7(33),\ldots,p_7(35)$, beyond the range $n\le32$ of GMMR's table \cite{GMMRData}.}\label{tab:p7new}\\
\toprule
$n$ & $p_7(n)$\\\midrule
\endfirsthead
\multicolumn{2}{@{}l}{\footnotesize Table~\thetable\ (continued)}\\
\toprule
$n$ & $p_7(n)$\\\midrule
\endhead
\midrule
\multicolumn{2}{r}{\footnotesize Continued on next page}\\
\endfoot
\bottomrule
\endlastfoot
33 & 35\,123\,727\,238\,126\,537\\
34 & 97\,715\,083\,517\,889\,540\\
35 & 270\,450\,168\,724\,504\,758\\
\end{longtable}
\endgroup

\subsection{Comparison with GMMR's posted data}\label{sec:gmmr-discrepancy}
Graffeo--Monavari--Moschetti--Ricolfi (GMMR) express the coefficients $p_m(n)$ through numerator polynomials $h_n(t)$ of their dimension-generating series, using their equation (1.3) \cite{GMMR2026}. Their arXiv paper determines these counts for $n\le30$ in every dimension, and its 2026 version on SSRN states them through $n=32$ \cite{GMMR2026}; their posted data file gives the polynomials $h_n(t)$ through $n=32$ \cite{GMMRData}. We compared that file with our coefficients for $m\le7$ and $n\le32$. The tables agree for every $n\le31$ and, at $n=32$, for $m\le4$. At $n=32$ our values are
\[
p_5(32)=21\,684\,370\,682\,500,\quad p_6(32)=637\,062\,384\,490\,959,\quad p_7(32)=12\,557\,013\,410\,901\,788.
\]
Write $\Delta_m=p_m^{\mathrm{GMMR}}(32)-p_m^{\mathrm{ours}}(32)$. We find
\[
\Delta_5=\delta,\qquad \Delta_6=7\delta,\qquad \Delta_7=28\delta,\qquad \delta=25\,422\,960.
\]
Thus, in the compared dimensions, $\Delta_m=\binom{m+1}{6}\delta$. In particular, the value of $p_7(32)$ obtained from GMMR's $h_{32}(t)$ is $12\,557\,014\,122\,744\,668$, which exceeds ours by $711\,842\,880$. Inverting the embedding-dimension transform, all counts $y_{32}^k$ agree for $0\le k\le8$ except $k=6$, where
\[
y_{32}^{6,\mathrm{ours}}=19\,142\,721\,504\,878,\qquad
y_{32}^{6,\mathrm{GMMR}}=19\,142\,746\,927\,838.
\]
Among embedding dimensions $0\le k\le8$ checked here, the discrepancy is confined to $k=6$. This comparison does not exclude further discrepancies at larger embedding dimensions. Our values at $n=32$ were each obtained under at least two conditioning decompositions (Section~\ref{sec:runs}), all implementing the same recurrence; this agreement does not independently settle which value is correct. A verification of $y_{32}^6$, equivalently $p_5(32)$, by a mathematically independent method would resolve the disagreement. Their data do not directly test the new coefficients above weight 32. The four authors were informed of the discrepancy by email on 11 November 2026, with the numerical comparison and a request to check the embedding-dimension-six count. Their data file as received and the comparison script accompany the computational certificate.

\subsection{How the coefficients were computed}\label{sec:runs}
All coefficients were computed in exact integer arithmetic, one CPU core per run. Permuting the $d$ coordinate axes of all slices together preserves weights and containment, so the contribution of a top block in \eqref{eq:top} is constant on its orbit under $S_d$; the runs therefore sum over orbit representatives weighted by orbit size, and memoise states under a canonical relabelling of axes. The ancillary program \texttt{hdfm.cpp} (Appendix~\ref{app:anc}) implements Algorithm D with first-member transfers for every chain, diagonal side states keyed by their residual floors (Appendix~\ref{app:residual}) and 128-bit coefficients, and it reproduces every table: on one 2.1\,GHz core it computes $p_3(0),\ldots,p_3(100)$ in about 3.6 hours and each of the tables in dimensions 4--7 in under 24 minutes. The solid-partition table was also computed by a separately written program specialised to solid partitions: once with first-member transfers, once evaluating the floor-free stationary chains by corner inclusion--exclusion, an exact identity outside the analysed algorithm, and twice with an earlier version of that program. All four runs agree with \texttt{hdfm.cpp} in every coefficient. In dimensions 4--7 the general-dimensional program from which \texttt{hdfm.cpp} is derived computed each table under two interior decompositions, using inclusion--exclusion for the floor-free stationary chains, again with identical results. For comparison, the computation of $p_3$ through weight 68 took about 360,000 CPU hours \cite{BGP}, and that through weight 72 about half a million \cite{DestainvilleGovindarajan}; these are historical resource reports, not a hardware-normalised comparison. All computed values agree with the published tables \cite{OEISSolid,OEIS4,OEIS5,OEIS6,GMMRData} wherever they overlap, apart from the discrepancy at $n=32$ described above, and with direct enumeration of order ideals at small weights. The runs use fixed conditioning depths and bounded caches with recomputation, so Theorems~\ref{thm:row-main} and \ref{thm:diag-main} do not bound their running times. Agreement between runs checks the implementations, not the recurrences, whose correctness is proved in Sections~\ref{sec:conditioning}--\ref{sec:diagonal}. The computational certificate among the ancillary files records every run, with its command, resources and checksums.

\appendix
\section{A complete hand-checkable solid-partition coefficient}\label{app:small}
This appendix supplies the full low-weight checks complementary to the hand calculation of Section~\ref{sec:example}. We first follow the complete decomposition in a case small enough to check directly. For a solid partition the alphabet is $\{0,1,2,3\}$. The valid words of lengths $0,1,\ldots,7$ number
\[
1,\ 4,\ 10,\ 26,\ 59,\ 140,\ 307,\ 684,
\]
namely $p_3(1),\ldots,p_3(8)$. For example, $321020$ has the three blocks $321$, $2$ and $\epsilon$: the last block represents a one-cell third slice, not an absent slice.

Cut a solid-partition word at every zero and condition on its third block $b$ ($K=3$). If the third block is present, the two blocks preceding it form a nested pair above $b$, and the suffix forms a chain inside $b$. A block of letter length $L(b)=m$ encodes a slice of weight $\omega(b)=m+1$; here and below $L(b)$ denotes the number of letters, whereas $\omega(b)$ denotes the weight of the encoded partition. Thus the contribution of a specified third block is
\begin{equation}\label{eq:word-example-decomp}
 q^{3\omega(b)} C(b;q)F(b;q),
\end{equation}
where $C(b;q)$ counts the two prefix slices by their total excess over two copies of $b$, and $F(b;q)$ counts all possible succeeding slices by their number of cells. Words with fewer than three blocks form a separate case.

For $n=7$, the third block can only be absent, empty ($\epsilon$), or one letter long. The contributions are as follows:
\begin{center}
\begin{tabular}{@{}p{36mm}p{78mm}r@{}}\toprule
Third block & Counting argument & Number\\\midrule
Absent & At most one zero: 86 words with no zero and 120 with one zero & 206\\
$\epsilon$ & Prefix coefficients $1,3,9,22,54$ at excess $0,\ldots,4$; suffix $1/(1-q)$ & 89\\
$1,2,3$ & Four words for each of three one-letter blocks & 12\\\midrule
Total & $206+89+12$ & \textbf{307}\\\bottomrule
\end{tabular}
\end{center}
For instance, the four words with third block $b=2$ are
\[
202020,\qquad210202,\qquad220202,\qquad320202.
\]
The first has three two-cell slices followed by a one-cell slice; the other three enlarge the first slice by one cell in one of three possible directions. Coordinate symmetry produces the corresponding lists for $b=1,3$.

\paragraph{The suffix.} Let $v$ be the slice encoded by $b$. If the suffix has first slice $u\subeq v$, it contributes $q^{\wt{u}}F(u;q)$. Hence
\begin{equation}\label{eq:word-example-suffix}
 F(v;q)=\sum_{u\subeq v}q^{\wt{u}}F(u;q),\qquad F(\varnothing;q)=1.
\end{equation}
For the one-cell slice $F(\epsilon;q)=(1-q)^{-1}$. For a two-cell slice, by symmetry,
\[
 F(1;q)=F(2;q)=F(3;q)=\frac1{(1-q)(1-q^2)}.
\]
Thus the words with $b=\epsilon$ are the coefficient of $q^4$ in $C(\epsilon;q)/(1-q)$.

\paragraph{The prefix as coupled diagonal words.} The two prefix slices are plane partitions. Each may be read along its main diagonal and its two sides; a diagonal is an ordinary partition and hence another word. In this sublanguage $\epsilon=(1)$, $1=(2)$, $0=(1,1)$, while $\varnothing$ means that a diagonal is absent. Write $\alpha_1\mid\alpha_2$ for a pair of diagonals of the outer and inner prefix slices. Their interlacing rules are local: the next pair lies between an explicit least successor and the current pair, subject to dominance of the outer member over the inner member. The centre supplies a minimum completion; the remaining cost splits between the positive and negative sides.

With third slice $\epsilon$, the two prefix slices must each contain the origin, and the eleven admissible central pairs at excess at most four are
\[
\begin{gathered}
(\epsilon\mid\epsilon),\ (1\mid\epsilon),\ (1\mid1),\ (11\mid\epsilon),\ (0\mid\epsilon),\\
(11\mid1),\ (111\mid\epsilon),\ (10\mid\epsilon),\ (11\mid11),\ (111\mid1),\ (1111\mid\epsilon).
\end{gathered}
\]
The positive and negative sides have the same unconstrained transfer series $T$ after the minimum has been accounted for. In particular, the side-transfer repeat accounts for the geometric factors in
\[
 T(\epsilon\mid\varnothing)=\frac1{1-q},\qquad
 T(\epsilon\mid\epsilon)=\frac1{(1-q)(1-q^2)},
\]
\[
 T(1\mid\varnothing)=\frac1{(1-q)(1-q^2)},\qquad
 T(1\mid\epsilon)=1+q+3q^2+4q^3+O(q^4).
\]
The last example follows by enumerating its four nonrepeating successors together with the repeat; the general handling of repeats is given in Appendix~\ref{app:residual}. Every one of the eleven possible centres at this budget is displayed below. A vector $a_0,a_1,\ldots$ denotes coefficients of the side series through the indicated residual budget; the final column is its square, because the two sides are identical.
\begin{center}\small
\begin{tabular}{@{}lcll@{}}\toprule
Centre $\alpha_1\mid\alpha_2$ & $\Delta_0$ & $T$ through degree $4-\Delta_0$ & $T^2$ through degree $4-\Delta_0$\\\midrule
$\epsilon\mid\epsilon$ & 0 & $1,1,2,2,3$ & $1,2,5,8,14$\\
$1\mid\epsilon$ & 1 & $1,1,3,4$ & $1,2,7,14$\\
$1\mid1$ & 2 & $1,1,3$ & $1,2,7$\\
$11\mid\epsilon$ & 2 & $1,1,3$ & $1,2,7$\\
$0\mid\epsilon$ & 3 & $1,2$ & $1,4$\\
$11\mid1$ & 3 & $1,1$ & $1,2$\\
$111\mid\epsilon$ & 3 & $1,1$ & $1,2$\\
$10\mid\epsilon$ & 4 & $1$ & $1$\\
$11\mid11$ & 4 & $1$ & $1$\\
$111\mid1$ & 4 & $1$ & $1$\\
$1111\mid\epsilon$ & 4 & $1$ & $1$\\\bottomrule
\end{tabular}
\end{center}
Summing $q^{\Delta_0}T^2$ over these centres gives
\begin{equation}\label{eq:prefix-small}
 C(\epsilon;q)=1+3q+9q^2+22q^3+54q^4+O(q^5).
\end{equation}
For a one-letter block $b$, the least prefix is two copies of $b$ and there are three excess-one enlargements, so
\[
 C(b;q)=1+3q+O(q^2).
\]
Equations \eqref{eq:word-example-suffix} and \eqref{eq:prefix-small} give all three parts of the computation.

\paragraph{A binding floor at weight eight.}
For $n=8$ the third block $b=2$ represents a horizontal domino. Its second cell lies on the first positive diagonal. For the centre $(\epsilon\mid\epsilon)$, the first positive-side transition is therefore forced by the floor to $(\epsilon\mid\epsilon)$ at zero excess cost, unlike the floor-free transition. In terms of a central entry $x$ at the affected position, its least-completion cost becomes
\[
 \Psi_1(x)=3x+(1-x)_+,
\]
which explicitly exhibits the binding floor. The excess-one prefixes of $b=2$ are $B_1=22,21,32$ with $B_2=2$, showing the three possible directions in which a cell can be added. Here the centre $(\epsilon\mid\epsilon)$ has central excess zero and the two side transfers agree with the floor-free $T(\epsilon\mid\epsilon)$, so their product contributes $1+2q+5q^2$. The centre $(1\mid\epsilon)$ has central excess one; its positive-side transfer is $T(\epsilon\mid\epsilon)+qT(1\mid\epsilon)=1+2q+O(q^2)$, while the negative-side transfer is $1+q+O(q^2)$. Its contribution is therefore $q+3q^2$. Three further central pairs have excess two and contribute $3q^2$ in total. This accounts explicitly for
\[
 C(2;q)=1+3q+11q^2+O(q^3).
\]
One of those last three pairs, $(0\mid\epsilon)$, has central excess two rather than the floor-free value three: the effective floor changes the increment to $\Psi_1(1)-\Psi_1(0)=2$. Since $F(2;q)=((1-q)(1-q^2))^{-1}=1+q+2q^2+O(q^3)$, it follows that
\[
 C(2;q)F(2;q)=1+4q+16q^2+O(q^3).
\]
By symmetry, $b=1,3$ give the same results. Thus the full count at $n=8$ is
\[
 p_3(8)=(160+267)+209+3\cdot16=427+209+48=684.
\]
The $n=7$ calculation demonstrates the exterior conditioning and diagonal-transfer machinery, and the $n=8$ case demonstrates a nontrivial floor. Intermediate-level family conditioning first occurs in dimensions four and higher, where the prefix words are recursively conditioned again until the diagonal stage is reached.

\section{Residual floors and level-independent diagonal states}\label{app:residual}
This appendix proves that the full plane-partition floor need not be carried at every step of the practical diagonal transfer. It is not needed for the theorems, but the program that computed the tables of Section~\ref{sec:computation} keys its diagonal side states by these residual floors, so their exactness rests on it.

Fix one plane-partition member $\phi$ of the floor, a side level $\ell\ge0$, and a current diagonal $\alpha\supeq d_\ell(\phi)$. Let
\[
 D_\ell=\{(a,b)\in\NN^2:b-a>\ell\}
\]
be the region strictly beyond the current positive diagonal, and define its \emph{residual floor steps}
\begin{equation}\label{eq:residual}
 R_\ell(\alpha)=\{(a,b)\in D_\ell:
 \phi(a,b)>\phi(a+1,b),\ \phi(a,b)>\alpha(b-\ell)\},
\end{equation}
carrying the numerical floor values $\phi(a,b)$. For transfer and memoisation we record a cell relative to the diagonal as $(s,a)$ with $s=b-a-\ell-1$. For a $\Gamma$-family, retain one such labelled residual set for each member; the coupling inequalities are already encoded in the family of current diagonals.

\begin{lemma}[A surviving step in each binding column]\label{lem:columnstep}
If $(a,b)\in D_\ell$ and $\phi(a,b)>\alpha(b-\ell)$, let $a^*$ be the largest integer in $[a,b-\ell-1]$ such that $\phi(a^*,b)=\phi(a,b)$. Then $(a^*,b)\in R_\ell(\alpha)$.
\end{lemma}
\begin{proof}
If $a^*+1\le b-\ell-1$, maximality gives $\phi(a^*,b)>\phi(a^*+1,b)$. If $a^*+1=b-\ell$, the latter entry lies on the current floor diagonal and hence is at most $\alpha(b-\ell)<\phi(a^*,b)$. In either case the vertical floor step is strict, and its value exceeds the current diagonal lower bound, so it lies in $R_\ell(\alpha)$.
\end{proof}

\begin{lemma}[Residual-floor sufficiency]\label{lem:residual}
Consider fillings $u$ of $D_\ell$ weakly decreasing along both coordinate directions and joining $\alpha$ according to the interlacing inequalities. Such a filling satisfies the full floor $u\supeq\phi$ on $D_\ell$ if and only if it satisfies $u(a,b)\ge\phi(a,b)$ for every $(a,b)\in R_\ell(\alpha)$. Consequently its continuation series depends only on $(\alpha,R_\ell(\alpha))$ in the relative coordinate system, independently of the absolute level and of the other values of $\phi$.
\end{lemma}
\begin{proof}
Necessity is immediate. For sufficiency, fix $(a,b)\in D_\ell$. Monotonicity down column $b$ and the junction with the current diagonal force $u(a,b)\ge\alpha(b-\ell)$. This meets the floor whenever $\phi(a,b)\le\alpha(b-\ell)$. Otherwise Lemma~\ref{lem:columnstep} produces a binding step $(a^*,b)\in R_\ell(\alpha)$ with $a^*\ge a$ and $\phi(a^*,b)=\phi(a,b)$. The filled array is weakly decreasing in its first coordinate, so
\[
 u(a,b)\ge u(a^*,b)\ge\phi(a^*,b)=\phi(a,b).
\]
Thus every omitted floor condition is redundant. In coordinates relative to $\ell$, the region, junction and coupling conditions are fixed by $\alpha$, and the only effective additional constraints are the labelled cells of $R_\ell(\alpha)$.
\end{proof}

\begin{proposition}[Transfer and update from residual data]\label{prop:residual-transfer}
Let $R=R_\ell(\alpha)$, expressed in relative coordinates $(s,a)$ with floor values.
\begin{enumerate}[label=(\alph*),leftmargin=2em]
\item The least successor satisfies
\[
 \beta_i=\max\{\alpha(i+1),v_i\},
\]
where $v_i$ is the value of the residual cell $(s,a)=(0,i)$, or zero if it is absent.
\item For each $x\ge\beta_i$, the correction in $f_i(x)$ from \eqref{eq:f-side} is determined by $R$. More precisely, let $\rho_i(a)$ be the largest residual value in the physical column $i+\ell+1$ among rows $a,\ldots,i-1$, or zero if none exists. Then
\[
 \sum_{a<i}(\phi(a,i+\ell+1)-x)_+
 =\sum_{a<i}(\rho_i(a)-x)_+.
\]
Thus the successor set and $\Delta$ depend only on $(\alpha,R)$.
\item For every admissible successor $\alpha'$, the new residual data consist of the old cells with $s\ge1$ and value exceeding $\alpha'(a+s)$, shifted from $s$ to $s-1$.
\end{enumerate}
\end{proposition}
\begin{proof}
For (a), the first off-centre cell at $(i,i+\ell+1)$ either has floor value at most $\alpha(i+1)$, in which case it does not affect $\beta_i$, or it is a strict binding step by Lemma~\ref{lem:columnstep}, in which case it appears at $(0,i)$ in $R$ with precisely the needed value.

For (b), monotonicity of $\phi$ down a column gives $\rho_i(a)\le\phi(a,i+\ell+1)$. If the latter does not exceed $x$, both positive parts vanish. Otherwise $\phi(a,i+\ell+1)>x\ge\beta_i\ge\alpha(i+1)$; Lemma~\ref{lem:columnstep} supplies a residual step $a^*$ at or below $a$, still strictly above the next diagonal. Because the floor value at row $i$ is at most $\beta_i\le x$, this step has $a^*<i$ and the same value as the original cell. Hence $\rho_i(a)$ equals the original floor value whenever the positive part is nonzero, proving the equality.

For (c), the new future region $D_{\ell+1}$ is precisely the old region with offset $s\ge1$. A cell in it is a strict vertical step of the unchanged floor and is binding at the next level exactly when its value exceeds the next diagonal entry $\alpha'(b-\ell-1)=\alpha'(a+s)$. Such a cell was already binding at level $\ell$ because $\alpha'(a+s)\ge\beta(a+s)\ge\alpha(a+s+1)$. Therefore no new residual step is created: retain only the indicated old cells and decrease their offset by one.
\end{proof}

\begin{remark}[Repeats]
With level-independent keys $(\alpha,R)$, a cycle can occur only when $R$ is empty and the diagonal repeats unchanged. All positive offsets in nonempty $R$ decrease at each transition, while a nonrepeat successor reduces the diagonal size. When $R=\varnothing$ the repeat has cost $\wt\alpha>0$, so its self-contribution may be summed geometrically:
\begin{equation}\label{eq:repeat}
 T(\alpha;q)=\frac{\displaystyle\sum_{\alpha'\ne\alpha}q^{\Delta(\alpha')}T(\alpha';q)}{1-q^{\wt\alpha}}.
\end{equation}
At $\alpha=0$ the transfer is $1$. This is an implementation identity, not an additional asymptotic estimate.
\end{remark}

\section{Ancillary files}\label{app:anc}
The ancillary files of this article contain the program, the tables and the record of the computations of Section~\ref{sec:computation}.
\begin{center}\small
\begin{tabularx}{\textwidth}{@{}l X@{}}\toprule
File & Contents\\\midrule
\texttt{hdfm.cpp} & C++17 source of Algorithm D with first-member transfers and 128-bit integer coefficients; accepts $3\le d\le7$ and $n\le159$ (1350 lines)\\
\texttt{partitions\_hdfm.ipynb} & a Jupyter notebook that writes, compiles and runs the same program for given $d$ and $n$\\
\texttt{p3\_0\_100.txt}, \ldots, \texttt{p7\_0\_35.txt} & the tables $p_3(0..100)$, $p_4(0..50)$, $p_5(0..45)$, $p_6(0..40)$ and $p_7(0..35)$, one line ``$n$ $p_d(n)$'' per weight\\
\texttt{computational\_certificate.pdf} & every run with its command, host, resources and checksums, and the comparisons between runs, with published values and with direct enumeration\\
\texttt{RWC\_computational\_certificate.zip} & the sources, logs and outputs of all the runs, and a script that rechecks every comparison in the certificate\\\bottomrule
\end{tabularx}
\end{center}
The stated accepted input range is not a certificate of overflow-free exact arithmetic or successful completion throughout that range; the computation records certify the coefficient ranges listed above. The source \texttt{hdfm.cpp} has SHA-256
\begin{center}\scriptsize\texttt{b82de5d4b5d6c523be5f1fbd31533d46386fc60621fa07f3afc3bc4ec8704cb7}\end{center}
and is built with \texttt{g++ -O2 -march=native -std=c++17 -o hdfm hdfm.cpp}. Each of the following commands prints one table, one line ``$n$ $p_d(n)$'' per weight; its output is byte-identical to the corresponding table file, whose SHA-256 is shown below the command.
\begin{itemize}[leftmargin=1.5em]\small
\item $p_3(0..100)$: \texttt{./hdfm 3 100 3 --sigorder --memcap 5500}\\
{\scriptsize\texttt{b17fb2bcda1fe7e863c1d7477b7796fd2034fb2106e9b22080d08526adffc4a7}}
\item $p_4(0..50)$: \texttt{./hdfm 4 50 2 3 --canon-desc}\\
{\scriptsize\texttt{fe5a4b0128e5d3a363e9986efae64b0a8705e6bc243f909bebbd45d8cf05fef6}}
\item $p_5(0..45)$: \texttt{./hdfm 5 45 2 3 2 --canon-desc}\\
{\scriptsize\texttt{1b13c69e0346c16ae4d06505bd283be0e6d3ce57536f114d2f8c92286d6b13df}}
\item $p_6(0..40)$: \texttt{./hdfm 6 40 2 3 2 2 --canon-desc}\\
{\scriptsize\texttt{4780d6d7f620f41739182990ca650e0c30c622ceec13d6b2398ee7d737d8e50f}}
\item $p_7(0..35)$: \texttt{./hdfm 7 35 2 2 2 2 2}\\
{\scriptsize\texttt{21035c0acb4b80925051691ad6f0fcb2011bea216b359137a9d3b995bb631235}}
\end{itemize}
On one 2.1\,GHz core the first command takes about 3.6 hours and 5.8\,GiB of memory; the option \texttt{--memcap} sets the cache size at which the program clears its caches and recomputes, trading time for memory without changing the result. The other commands take under 24 minutes and at most 1.9\,GiB each.

\section*{Acknowledgments and research-method disclosure}
\addcontentsline{toc}{section}{Acknowledgments and research-method disclosure}
The initial idea of interior word conditioning led to a series of trial-and-error computational experiments. Algorithm design, C++ and Python implementation and optimisation, derivation and checking of the recursive identities, proof development, proof audits and independent software verification were conducted through sustained substantive interactions with OpenAI Astra (Medium), OpenAI Sol 6.1 and Claude Opus 5.5. The named AI systems are acknowledged as research assistants; responsibility for the paper, proofs, reproducibility record and any errors remains with the author.

\end{document}